\documentclass{amsart}

\newtheorem{theorem}{Theorem}[section]

\theoremstyle{definition}

\newtheorem{example}[theorem]{Example}

\theoremstyle{remark}

\begin{document}

\title[Logarithmic integrals]{A class of logarithmic integrals}

\author{Luis A. Medina}
\address{Department of Mathematics,
Tulane University, New Orleans, LA 70118}
\email{lmedina@math.tulane.edu}

\author{Victor H. Moll}
\address{Department of Mathematics,
Tulane University, New Orleans, LA 70118}
\email{vhm@math.tulane.edu}


\date{\today}


\begin{abstract}
We present a systematic study of integrals of the form 
\begin{equation}
I_{Q} = \int_{0}^{1}  Q(x) \log \log \frac{1}{x}  \, dx, 
\nonumber
\end{equation}
\noindent
where $Q$ is a rational function.  
\end{abstract}

\maketitle

\newcommand{\nn}{\nonumber}
\newcommand{\ba}{\begin{eqnarray}}
\newcommand{\ea}{\end{eqnarray}}
\newcommand{\ift}{\int_{0}^{\infty}}
\newcommand{\ifft}{\int_{- \infty}^{\infty}}
\newcommand{\ione}{\int_{0}^{1}}
\newcommand{\lx}{\log \log 1/x}
\newcommand{\al}{\omega}
\newcommand{\eps}{\epsilon}
\newcommand{\no}{\noindent}
\newcommand{\realpart}{\mathop{\rm Re}\nolimits}
\newcommand{\imagpart}{\mathop{\rm Im}\nolimits}
\newcommand{\st}{{ \quad \quad \quad  }}

\newtheorem{Definition}{\bf Definition}[section]
\newtheorem{Thm}[Definition]{\bf Theorem}
\newtheorem{Example}[Definition]{\bf Example}
\newtheorem{Lem}[Definition]{\bf Lemma}
\newtheorem{Note}[Definition]{\bf Note}
\newtheorem{Cor}[Definition]{\bf Corollary}
\newtheorem{Prop}[Definition]{\bf Proposition}
\newtheorem{Problem}[Definition]{\bf Problem}
\numberwithin{equation}{section}

\section{Introduction} \label{sec-intro}
\setcounter{equation}{0}

The classical table of integrals by I. S.
Gradshteyn and I. M. Ryzhik \cite{gr} contains
very few evaluations of the form
\begin{equation}
I_{Q} := \ione Q(x) \lx \, dx, 
\label{type-1}
\end{equation}
\noindent
where $Q$ is a rational function. The example 
\begin{equation}
\ione \frac{1}{x^{2}+1} \lx  \, dx = 
\frac{\pi}{2} \log \left( \frac{\sqrt{2 \pi} \, \Gamma \left( \tfrac{3}{4} \right)}
{ \Gamma \left( \tfrac{1}{4} \right)} \right),
\label{int-vardi}
\st
\end{equation}
\noindent
written in its trigonometric version
\begin{equation}
\ione \frac{1}{x^{2}+1} \lx  \, dx = 
\int_{\pi/4}^{\pi/2} \log \log \tan x \, dx,
\label{trig-version}
\st
\end{equation}
\noindent
is the subject of Vardi's remarkable paper \cite{vardi1}. This example
appears as $4.229.7$ in \cite{gr}
and also in the equivalent form 
\begin{equation}
\ift \frac{\log x \, dx }{\cosh x} = 
\pi \log \left( \frac{\sqrt{2 \pi}  \, \Gamma \left( \tfrac{3}{4} \right)}
{ \Gamma \left( \tfrac{1}{4} \right)} \right),
\st
\end{equation}
\noindent
as $4.371.1$ in the same table. 

We present here a systematic study of the logarithmic integrals (\ref{type-1}).
Throughout the paper we indicate whether Mathematica 6.0 is capable of 
evaluating the integrals considered. For example, a direct symbolic evaluation 
gives (\ref{int-vardi}) as
\begin{equation}
\ione \frac{1}{x^{2}+1} \lx  \, dx = 
\frac{\pi}{4} \log \left( \frac{ 4 \pi^{3}}{\Gamma( \tfrac{1}{4} )^{4}} 
\right). 
\st
\end{equation}
\noindent
The reader should be aware that the question of whether a definite 
integral is computable by a symbolic language depends on the form
in which the integrand is expressed. For
instance, Mathematica 6.0 is unable to evaluate the 
trigonometric version of (\ref{int-vardi}) given as the right-hand side of 
(\ref{trig-version}). 

The idea exploited here, introduced by I. Vardi in \cite{vardi1},  is to 
associate to each function $Q$ a gamma factor
\begin{equation}
\Gamma_{Q}(s) := \ione Q(x) \left( \log \frac{1}{x} \right)^{s-1} \, dx,
\end{equation}
\noindent
so that, the integral (\ref{type-1}) is given by 
\begin{equation}
I_{Q} = \Gamma_{Q}'(1).
\end{equation}

An explicit evaluation of $I_{Q}$ is achieved in the case where $Q(x)$ is
analytic at $x=0$. Starting with the expansion 
\begin{equation}
Q(x) = \sum_{n=0}^{\infty} a_{n} x^{n},
\end{equation}
\noindent
we associate an $L$-series 
\begin{equation}
L_{Q}(s) := \sum_{n=0}^{\infty} \frac{a_{n}}{(n+1)^{s}}.
\end{equation}
\noindent
The integral (\ref{type-1}) is now evaluated as
\begin{equation}
I_{Q} = - \gamma L_{Q}(1) + L_{Q}'(1),
\label{vardi-key}
\end{equation}
\noindent
where $\gamma$ is the {\em Euler-Mascheroni} constant. The
interesting story of this fundamental constant can be found
in \cite{havil1}.  The identity (\ref{vardi-key})  is essentially Vardi's 
method for the evaluation of (\ref{int-vardi}). Naturally, to obtain an
{\em explicit} evaluation of $I_{Q}$, one needs to express $L_{Q}(1)$ and 
$L_{Q}'(1)$ in terms of special functions.  We employ here the Riemann 
zeta function
\begin{equation}
\zeta(s) = \sum_{n=1}^{\infty} \frac{1}{n^{s}}, \quad s > 1,
\end{equation}
\noindent
and its alternating form
\begin{equation}
\zeta_{a}(s) = \sum_{n=1}^{\infty} \frac{(-1)^{n}}{n^{s}} = 
-(1 - 2^{1-s}) \zeta(s), \quad s > 0.
\st
\end{equation}
\noindent
The relations 
\begin{eqnarray}
\zeta(2n) & = & \frac{(-1)^{n+1} ( 2 \pi)^{2n} B_{2n}}{2(2n)!}, \quad
n \in {\mathbb{N}}_{0} := \mathbb{N} \cup \{ 0 \},  \st \label{zeta-even} \\
\zeta(1-n) & = & \frac{(-1)^{n+1} B_{n}}{n}, \quad 
n \in \mathbb{N}, \st \label{zeta-neg} \\
\zeta'(-2n) & = & (-1)^{n} \frac{(2n)! \, \zeta(2n+1) }{2 (2 \pi)^{2n}}, 
\quad n \in \mathbb{N}, \st \label{zetaprime} \\
\zeta'(0) & = & - \log \sqrt{2 \pi}, \st \label{zeta-0}
\end{eqnarray}
\noindent
where $B_{n}$ are the Bernoulli numbers, will be 
used to simplify the integrals discussed below. 

A second function that is used in the evaluations 
described here is the {\em polylogarithm function} defined by 
\begin{equation}
\text{PolyLog}[c,x] := \sum_{n=1}^{\infty} \frac{x^{n}}{n^{c}}
\label{polylog-def0}
\end{equation}
\noindent
and its derivative
\begin{eqnarray}
\text{PolyLog}^{(1,0)}[c,x] & := & \frac{d}{dc} \text{PolyLog}[c,x]  
\label{polyder-def0} \\
& = & - \sum_{n=1}^{\infty} \frac{\log n}{n^{c}}x^{n}.
\nonumber
\end{eqnarray}

\medskip

A complete description of $I_{Q}$ is determined here in the case
where $Q(x)$ is a rational function. The 
method of partial fractions shows that it suffices to consider three types of 
integrals: \\

\noindent
the first type is 
\begin{equation}
P_{j} := \ione x^{j} \lx \, dx, \label{Pm}
\end{equation}
\noindent
that gives the polynomial part of $Q$,  \\

\noindent
the second type is 
\begin{equation}
R_{m,j}(a) := \ione \frac{x^{j} \, \lx }{(x+a)^{m+1}}
 \, dx, \label{Bmj}
\end{equation}
\noindent
that treats the real poles of $Q$, and \\

\noindent
the third type is 
\begin{equation}
C_{m,j}(a,b) := \ione \frac{x^{j} \, \lx }{(x^{2} + ax+b)^{m+1}}
\, dx, \label{Cmj}
\end{equation}
\noindent
with $a^{2} - 4b<0$. This last case deals with the non-real poles of $Q$. 

\medskip

\noindent
{\em Integrals of first type}. These are simple. They are 
evaluated in (\ref{value-Pm}) as 
\begin{equation}
P_{j} := \ione x^{j} \lx  \, dx = -\frac{\gamma+\log(j+1)}{j+1}. 
\nonumber
\st
\end{equation}

\medskip

\noindent
{\em Integrals of second type}. 
The special case $R_{m,0}(1)$ is evaluated first. We introduce the 
polynomial
\begin{equation}
T_{m}(x) := \sum_{j=0}^{m} (-1)^{j} A_{m+1,j+1}x^{j},
\end{equation}
\noindent
with $A_{m,j}$ the Eulerian numbers given in (\ref{eulerian-1}). Then 
the integral
\begin{equation}
E_{m} := \ione \frac{T_{m-1}(x) \, \lx }{(x+1)^{m+1}} \, dx
\end{equation}
\noindent
is evaluated, for $m > 1$,  as 
\begin{equation}
E_{m} = (1-2^{m})\zeta'(1-m) +(-1)^{m} \left( \gamma ( 2^{m}-1) + 2^{m} \log 2 \right)
\frac{B_{m}}{m}.
\label{expr-em0}
\st
\end{equation}
\noindent
The sequence $E_{m}$ is then used to 
produce a recurrence for $R_{m,0}(1)$. The initial condition 
\begin{equation}
R_{0,0}(1) = - \frac{1}{2} \log^{2}2
\st
\end{equation}
\noindent
is given in Example \ref{r00}. Then, the values of $R_{m,0}(1)$ are
obtained from
\begin{equation}
R_{m,0}(1)  =  \frac{E_{m}}{b_{0}(m)} - 
\sum_{k=1}^{m-1} \frac{b_{k}(m)}{b_{0}(m)} R_{m-k,0}(1),
\end{equation}
\noindent
where 
\begin{equation}
b_{k}(m) = (-1)^{k} \sum_{j=0}^{m-1} \binom{j}{k} A_{m,j+1}.
\end{equation}
\noindent
This is described in  Corollary \ref{cor-recu}. 

The evaluation of $R_{m,0}(a)$, for $a \neq 1$, appears in Proposition 
\ref{prop-51} and Corollary \ref{cor-53}: the value of $R_{0,0}(a)$ is 
given by
\begin{equation}
R_{0,0}(a) := \ione \frac{\lx \, dx}{x+a} = - \gamma \log (1 + 1/a) - 
\text{PolyLog}^{(1,0)}[1,-1/a],
\end{equation}
\noindent
and, for $m > 0$, we have 
\begin{eqnarray}
R_{m,0}(a) & = & - \frac{\gamma}{a^{m}(1+a)m} - 
\frac{\gamma}{a^{m+1} m!} \sum_{j=2}^{m} \frac{S_{1}(m,j) T_{j-2}(1/a)}
{(1+ 1/a)^{j}} \nonumber \\
& & - \frac{1}{a^{m} m!} 
\sum_{j=1}^{m} S_{1}(m,j) \text{PolyLog}^{(1,0)}[1-j,-1/a].
\nonumber
\end{eqnarray}
\noindent
Here $S_{1}(m,j)$ are the (signless) Stirling numbers of the first kind 
defined by the expansion 
\begin{equation}
(t)_{m} = \sum_{j=1}^{m} S_{1}(m,j) t^{j},
\label{stir-1a}
\end{equation}
\noindent
where $(t)_{m} = t(t+1)(t+2) \cdots (t+m-1)$ is the Pochhammer symbol.
The function 
$\text{PolyLog}^{(1,0)}[c,x]$ is defined in (\ref{polyder-def0}).

\medskip

For $j>0$, the value of $R_{m,j}(a)$ is now obtained from the recursion in
Theorem \ref{bmj-rec}, written here as 
\begin{equation}
R_{m,0}(a) = \sum_{j=0}^{r} \alpha_{j,r}(a) R_{m-r+j,j}(a),
\end{equation}
\noindent
where 
\begin{equation}
\alpha_{j,r}(a) := (-1)^{j} \binom{r}{j} a^{-r}.
\end{equation}
\noindent
This can be used for increasing values of the free parameter $r$, to obtain 
analytic expressions for $R_{m,j}(a)$. For instance, $r=1$ gives 
\begin{equation}
R_{m,0}(a) = \alpha_{0,1}(a)R_{m-1,0}(a) + \alpha_{1,1}(a) R_{m,1}(a), 
\label{form-131}
\end{equation}
\noindent
that determines $R_{m,1}(a)$ in terms of $R_{m,0}(a)$ and $R_{m-1,0}(a)$, that
were previously computed.  The 
value $r=2$ gives 
\begin{equation}
R_{m,0}(a) = \alpha_{0,2}(a) R_{m-2,0}(a) + \alpha_{1,2}(a)R_{m-1,1}(a)  +
\alpha_{2,2}(a)R_{m,2}(a),
\label{form-132}
\end{equation}
\noindent
that determines $R_{m,2}(a)$ in terms of previously computed integrals.  This 
procedure determines all the integrals $R_{m,j}(a)$.  \\

\noindent
{\em Integrals of third type}. These are integrals where the corresponding 
quadratic factor has non-real zeros. The expression
\begin{equation}
x^{2} + ax + b = (x-c)(x- \bar{c}) = x^{2} - 2rx \cos \theta  + r^{2}
\end{equation}
\noindent
is used to define
\begin{equation}
\label{def-D0}
D_{m,j}(r, \theta) = \ione \frac{x^{j} \, \lx }
{(x^{2}-2rx \cos \theta  + r^{2})^{m+1}}
\, dx. 
\end{equation}
\noindent
Naturally $C_{m,j}(a,b) = D_{m,j}(r, \theta)$, we are simply emphasizing the 
polar representation of the poles.  \\

The value $D_{0,0}(r, \theta)$ is computed 
first. Theorem \ref{thm-73} treats the case $r=1$, with the value 
\begin{equation}
D_{0,0}(1, \theta ) 
 =  \frac{\pi}{2 \sin \theta} 
\left[ ( 1 - \theta/\pi) \log 2 \pi + 
\log \left( \frac{\Gamma(1- \theta/ 2 \pi)}{\Gamma( \theta/ 2 \pi)} 
\right) \right].
\nonumber 
\end{equation}
\noindent
The case $r \neq 1$ is given in Theorem \ref{thm-78} as 
\begin{eqnarray}
D_{0,0}(r, \theta) 
& = & - \frac{\gamma}{r \sin \theta} 
\tan^{-1} \left( \frac{\sin \theta}{r - \cos \theta}  \right) \nonumber \\
 &  & + \frac{1}{2 r i \sin \theta} 
\left( 
\text{PolyLog}^{(1,0)}[1, e^{i \theta}/r] - 
\text{PolyLog}^{(1,0)}[1, e^{-i \theta}/r] 
\right). 
\nonumber 
\end{eqnarray}

The next step is to compute $D_{0,1}(r, \theta)$. This is described at the 
end of Section \ref{sec-complex}. The result is expressed in terms of the 
{\em Lerch zeta function}
\begin{equation}
\Phi(z,s,a) := \sum_{k=0}^{\infty} \frac{z^{k}}{(k+a)^{s}}, 
\label{lerch-def}
\end{equation}
as 
\begin{eqnarray}
D_{0,1}(r, \theta) & = & - \frac{\gamma}{2} \log \left( \frac{r^{2} - 
2 r \cos \theta +1}{r^{2}} \right) \nonumber \\
& - & \gamma \cot \theta \, \tan^{-1} \left( \frac{\sin \theta}
{r - \cos \theta } \right) \nonumber \\
& + & 
\frac{1}{2ri \sin \theta } 
\left[ \Phi^{(0,1,0)} \left( \frac{e^{i \theta}}{r}, 1 , 1 \right) - 
       \Phi^{(0,1,0)} \left( \frac{e^{- i \theta}}{r}, 1 , 1 \right) \right]. 
\nonumber
\end{eqnarray}
\noindent
The reader will find in \cite{lerch1} information about this function. 

\medskip

The values of $D_{m,j}(r, \theta)$  for $m, \, j > 0$, 
are determined by the recurrences 
\begin{equation}
D_{m,j}(r, \theta)  =  - \frac{1}{2rm \sin \theta} 
\frac{\partial}{\partial \theta}
D_{m-1,j-1}(r, \theta), 
\label{recu-one}
\end{equation}
\noindent
and
\begin{equation}
D_{m,j}(r, \theta) = \frac{1}{2m \cos \theta} \left( 
\frac{\partial}{\partial r} D_{m-1,j-1}(r, \theta) + 2rm D_{m,j-1}(r, \theta)
\right). 
\label{recu-two}
\end{equation}
\noindent
These follow directly from the definition of $D_{m,j}(r, \theta)$. Details 
are given in Section \ref{sec-complex}. \\

\noindent
{\bf Comment}. Integration by parts, shows that the integrals $I_{Q}$ in 
(\ref{type-1}) include those of the form
\begin{equation}
J_{Q} = \ione \frac{Q(x) \, dx}{\log x}. 
\label{type-2}
\end{equation}
\noindent
This class was originally studied 
by V. Adamchik \cite{adamchik3}. They were considered by Baxter, Tempereley 
and Ashley \cite{baxter1} in their work on the so-called Potts model for the 
triangular lattice. In that model, the generating function has the form
\begin{equation}
P_{3}(t) := 3 \ift \frac{\sinh( (\pi-t)x) \, \sinh \left( \frac{2 t x}{3} 
\right) }{x \sinh( \pi x) \cosh( tx)} \, dx. 
\end{equation}
\noindent
V. Adamchik determined analytic expressions for (\ref{type-2}) in the case 
where the denominator of $Q$ is a cyclotomic polynomial. The expressions
involve derivatives of the Hurwitz zeta function 
\begin{equation}
\zeta(z,q) := \sum_{n=0}^{\infty} \frac{1}{(n+q)^{z}}. 
\label{hurwitz}
\end{equation}
\noindent
Zhang Nan-Hue and  K. Williams \cite{yue3}, \cite{yue2} used the 
Hurwitz zeta function to 
evaluate definite integrals similar to the class considered here. Examples
of integrals that involve $\zeta(z,q)$ in the integrand are given in 
\cite{espmoll1}, \cite{espmoll2} and \cite{espmoll3}.

\section{The main tool} \label{sec-main}
\setcounter{equation}{0}

The evaluation of the integral
\begin{equation}
I_{Q} := \ione Q(x) \lx  \, dx,
\label{type-1a}
\end{equation}
\noindent
for a general function $Q(x)$, 
is achieved by introducing the function
\begin{equation}
\Gamma_{Q}(s) := \ione Q(x) \left( \log 1/x \right)^{s-1} \, dx. 
\end{equation}
\noindent
The next result is elementary. 

\begin{Lem}
\label{eval-int}
The integral $I_{Q}$ is given by 
\begin{equation}
I_{Q} = \Gamma_{Q}'(1).
\end{equation}
\end{Lem}

\begin{example}
The simplest case is $Q(x) \equiv 1$. Here we obtain 
\begin{equation}
\ione \lx  \, dx = \Gamma'(1),
\st
\end{equation}
\noindent
where 
\begin{equation}
\Gamma(s) = \ione \left( \log \frac{1}{x} \right)^{s-1} \, dx = 
\ift t^{s-1}e^{-t} \, dt,
\st
\end{equation}
\noindent
is the classical {\em gamma function}. The reader will find in \cite{irrbook}
the identity
\begin{equation}
\Gamma'(1) = - \gamma,
\st
\end{equation}
\noindent
where $\gamma$ is the {\em Euler-Mascheroni} constant. This 
example appears as $4.229.1$ in \cite{gr}. 
\end{example}

\begin{example}
Consider now the case $Q(x) = x^{a}$, for $a \in \mathbb{R}$. Observe that
\begin{equation}
\Gamma_{x^{a}}(s) = \ift e^{-(a+1)t} t^{s-1} \, dt = 
\frac{\Gamma(s)}{(a+1)^{s}}, \text{ for } a > -1 \text{ and } s > 0. 
\label{form-1}
\st
\end{equation}
\noindent
Differentiate with respect to $s$ at $s=1$ to produce
\begin{equation}
\ione x^{a} \lx  \, dx = -\frac{\gamma+\log(a+1)}{a+1}.
\label{form-3}
\st
\end{equation}
\noindent
Differentiating (\ref{form-1}) $n$ times with respect to $s$, yields 
\begin{equation}
\ione x^{a} \log^{n}x \lx  \, dx = 
\frac{(-1)^{n+1} \, n!}{(1+a)^{n+1}} \left( \log(1+a) + \gamma -
H_{n} \right),
\label{form-2}
\st
\end{equation}
\noindent
where $H_{n}$ is the $n$-th harmonic number. Mathematica 6.0 is unable to
evaluate (\ref{form-1}) if both $a$ and $n$ are entered as  
parameters. The same holds for 
(\ref{form-2}). 
\end{example}

\begin{Note}
The expression (\ref{form-3}), with  $a= m \in \mathbb{N}$, provides
the evaluation of the integral $P_{m}$ in (\ref{Pm}): 
\begin{equation}
P_{m} := \ione x^{m} \lx  \, dx = -\frac{\gamma+\log(m+1)}{m+1}.
\label{value-Pm}
\st
\end{equation}
\noindent
This appears as $4.325.8$ in \cite{gr}. 
\end{Note}

\section{The case where $Q$ is analytic at $x=0$} \label{sec-analy}
\setcounter{equation}{0}

In this section we consider the evaluation of the integral
\begin{equation}
I_{Q} := \ione Q(x) \lx  \, dx, 
\label{type-1b}
\end{equation}
\noindent
where $Q$ is admits an expansion 
\begin{equation}
Q(x) = \sum_{n=0}^{\infty} a_{n}x^{n}.
\label{expansion}
\end{equation}
\noindent
The expression for $I_{Q}$ is expressed in terms of the associated 
$L$-function defined by
\begin{equation}
L_{Q}(s) := \sum_{n=0}^{\infty} \frac{a_{n}}{(n+1)^{s}}.
\label{L-def}
\end{equation}

The idea for the next lemma comes from \cite{vardi1}.

\begin{Lem}
The function $\Gamma_{Q}$ satisfies $\Gamma_{Q}(s) = \Gamma(s) L_{Q}(s)$.
\end{Lem}
\begin{proof}
The linearity of $\Gamma_{Q}(s)$ in the $Q$-argument shows that
\begin{equation}
\Gamma_{Q}(s) = \sum_{n=0}^{\infty} a_{n} \Gamma_{x^{n}}(s).
\end{equation}
\noindent
The result now follows from the value of $\Gamma_{x^{n}}(s)$ in (\ref{form-1}).
\end{proof}

\begin{Thm}
\label{main-thm}
Assume $Q$ is given by (\ref{expansion}). Then 
\begin{equation}
I_{Q} := \ione Q(x) \lx  \, dx = -\gamma L_{Q}(1) + 
L_{Q}'(1).
\label{main-formula}
\end{equation}
\end{Thm}
\begin{proof}
Differentiate the expression for $\Gamma_{Q}$ in the previous lemma and use
the result of Lemma \ref{eval-int}.
\end{proof}

The theorem reduces the evaluation of $I_{Q}$ to the evaluation of 
$L_{Q}(1)$ and $L'_{Q}(1)$. The first series of examples come from
prescribing the coefficients $a_{n}$ of $Q(x)$ so that the $L_{Q}$ function 
is relatively simple.

\begin{example}
\label{first}
Choose $a_{n} = 1/(n+1)$. Then 
\begin{equation}
Q(x) = \sum_{n=0}^{\infty} \frac{x^{n}}{n+1} = - \frac{\log(1-x)}{x}.
\st
\end{equation}
\noindent
Then
\begin{equation}
L_{Q}(s) = \sum_{n=0}^{\infty} \frac{1}{(n+1)^{s+1}} = \zeta(s+1),
\st
\end{equation}
\noindent
where $\zeta(s)$ is the classical Riemann zeta function. 
Theorem \ref{main-thm} gives
\begin{equation}
\ione \frac{\log(1-x)}{x} \lx  \, dx = 
\frac{\gamma \pi^{2}}{6} - \zeta'(2).
\label{exam-1}
\st
\end{equation}

Mathematica 6.0 is unable to check this directly, but the change of variables
$x = e^{-t}$ transforms (\ref{exam-1}) to 
\begin{equation}
\ift \log t \log(1 - e^{-t}) \, dt = 
\frac{\gamma \pi^{2}}{6} - \zeta'(2).
\st
\end{equation}
\noindent
This is computable by Mathematica 6.0.  \\

The constant 
\begin{equation}
\zeta'(2) = - \sum_{n=1}^{\infty} \frac{\log n}{n^{2}}
\st
\end{equation}
\noindent
can be expressed in terms of the {\em Glaisher constant}
\begin{equation}
\log A := \frac{1}{12} - \zeta'(-1),
\st
\label{glaisher}
\end{equation}
\noindent
by 
\begin{equation}
\zeta'(2) = \frac{\pi^{2}}{6} \left( \gamma + \log(2 \pi) - 12 \log A \right). 
\st
\end{equation}
\noindent
This gives
\begin{equation}
\ione \frac{\log(1-x)}{x} \lx  \, dx = 
\frac{\pi^{2}}{6} \left( 12 \log A - \log 2 \pi \right)
\label{exam-2}
\st
\end{equation}
\noindent
as an alternative form for (\ref{exam-1}).
\end{example}

\medskip

\begin{example}
\label{second}
We now consider the alternating version of Example \ref{first} and choose
$a_{n} = (-1)^{n}/(n+1)$. In this case 
\begin{equation}
Q(x) = \sum_{n=0}^{\infty} \frac{(-1)^{n} x^{n}}{n+1} = \frac{\log(1+x)}{x},
\st
\end{equation}
\noindent
and the corresponding $L$-series is
\begin{equation}
L_{Q}(s) = \sum_{n=0}^{\infty} \frac{(-1)^{n}}{(n+1)^{s+1}} =
(1- 2^{-s}) \zeta(s+1).
\st
\end{equation}
\noindent
Theorem \ref{main-thm} and the evaluations 
\begin{equation}
L_{Q}(1) = \frac{\pi^{2}}{12} \text{ and } 
L_{Q}'(1) = \frac{\pi^{2} \log 2}{12} + \frac{1}{2} \zeta'(2),
\st
\end{equation}
\noindent
now yield
\begin{equation}
\ione \frac{\log(1+x)}{x} \lx \, dx = 
\frac{\pi^{2}}{12} \left( \log 2 - \gamma \right) + \frac{1}{2} \zeta'(2).
\st
\end{equation}
\noindent
This can also be expressed as 
\begin{equation}
\ione \frac{\log(1+x)}{x} \lx \, dx = 
\frac{\pi^{2}}{12} \left( \log 4 \pi - 12 \log A \right). 
\st
\end{equation}
\noindent
As in the previous example, Mathematica 6.0 is unable to produce this 
evaluation, but it succeeds with the alternate version
\begin{equation}
\ift \log t \log(1 + e^{-t}) \, dt = 
\frac{\pi^{2}}{12} \left( \log 2 - \gamma \right)  + \frac{1}{2} 
\zeta'(2).
\st
\end{equation}
\end{example}

\medskip

\begin{example}
\label{third}
Adding the results of the first two examples yields
\begin{equation}
\ione \frac{\log(1-x^2)}{x} \lx \, dx = 
\frac{\pi^{2}}{12} \left( \log 2 + \gamma \right) -
\frac{1}{2} \zeta'(2).
\st
\end{equation}
\noindent
Their difference produces
\begin{equation}
\ift \log t \log \tanh t \, dt = \frac{\gamma \pi^{2}}{8} - 
\frac{3}{4} \zeta'(2) + \frac{\pi^{2} \log 2 }{12}. 
\st
\end{equation}
\noindent
This cannot be evaluated symbolically. 
\end{example}

\medskip

\begin{example}
\label{fourth}
This example generalizes Example \ref{first}. The integrand involves the 
{\em polylogarithm} function defined in (\ref{polylog-def0}). The 
choice $a_{n} = 1/(n+1)^{c}$ produces the function
\begin{equation}
Q(x) = \sum_{n=0}^{\infty} \frac{x^{n}}{(n+1)^{c}} = 
\frac{1}{x} \text{PolyLog}[c,x],
\st
\end{equation}
\noindent
and the corresponding $L$-function is 
\begin{equation}
L_{Q}(s) = \sum_{n=0}^{\infty} \frac{1}{(n+1)^{s+c}} = \zeta(s+c).
\st
\end{equation}
\noindent
Then (\ref{main-formula}) gives
\begin{equation}
\ione \frac{\text{PolyLog}[c,x] }{x} \lx \, dx = -\gamma \zeta(c+1) +
\zeta'(c+1).
\st
\end{equation}
\end{example}

\begin{example}
Choosing now $a_{n} = (-1)^{n}/(n+1)^{c}$ gives 
\begin{equation}
Q(x) = \sum_{n=0}^{\infty} \frac{(-1)^{n} x^{n}}{(n+1)^{c}} = 
- \frac{1}{x} \text{PolyLog}[c,-x]
\st
\end{equation}
\noindent
and 
\begin{equation}
L_{Q}(s) = \sum_{n=0}^{\infty} \frac{(-1)^{n}}{(n+1)^{s+c}} = 
(1 - 2^{1-s-c}) \zeta(s+c).
\st
\end{equation}
\noindent
We conclude that 
\begin{multline}
\ione \frac{\text{PolyLog}[c,-x] }{x} \lx \, dx  =  
\left( \gamma (1 - 2^{-c}) - 2^{-c} \log 2 \right) \zeta(c+1) \nonumber \\
  -  (1- 2^{-c}) \zeta'(c+1).  \st \nonumber
\end{multline}
\noindent
In the special case where $c$ is a negative integer, the function $Q$ reduces 
to a rational function. Details are provided in section \ref{sec-special}.
\end{example}

\medskip

\begin{example}
Interesting integrands can be produced when  $a_{n}$ that are 
periodic sequences. For example, the choice 
\begin{equation}
a_{n} = \frac{1}{n} \cos \left( \frac{2 \pi n}{3} \right) 
\end{equation}
\noindent
gives the function 
\begin{equation}
Q(x) = \sum_{n=1}^{\infty} 
\frac{1}{n} \cos \left( \frac{2 \pi n}{3} \right) x^{n-1} = 
- \frac{\log( 1 + x + x^{2})}{2x}, 
\st
\end{equation}
\noindent
and the corresponding $L$-function 
\begin{equation}
L_{Q}(s) = \sum_{n=1}^{\infty} \frac{\cos( \tfrac{2 \pi n}{3}) }{n^{s+1}} = 
- \frac{1-3^{-s}}{2} \zeta(s+1). 
\st
\end{equation}
\noindent
Then (\ref{main-formula}) gives
\begin{equation}
\ione \frac{\log(1 + x + x^{2})}{x} \lx \, dx = 
- \frac{\gamma \pi^{2}}{9} + \frac{1}{18} \pi^{2} \log 3 + 
\frac{2}{3} \zeta'(2). 
\st
\end{equation}
\end{example}

\begin{example}
This example presents a second periodic sequence. The choice  
\begin{equation}
a_{n} = \frac{1}{n} \cos \left( \frac{2 \pi n}{5} \right) 
\end{equation}
\noindent
gives the function 
\begin{equation}
Q(x) = \sum_{n=1}^{\infty} 
\frac{1}{n} \cos \left( \frac{2 \pi n}{5} \right) x^{n-1} = 
- \frac{\log( 1  - \varphi  x + x^{2})}{2x}, 
\st
\end{equation}
\noindent
where $\varphi = (\sqrt{5}-1)/2$. The corresponding $L$-function is 
\begin{equation}
L_{Q}(s) = \sum_{n=1}^{\infty} \frac{\cos( \tfrac{2 \pi n}{5}) }{n^{s+1}} = 
\frac{1}{2} \left(
\text{PolyLog}[s+1, e^{-2 \pi i/5}] + 
\text{PolyLog}[s+1, e^{2 \pi i/5}]
\right).
\nonumber 
\end{equation}
\noindent
The values $L_{Q}(1) = \pi^{2}/150$ and 
\begin{equation}
L_{Q}'(1) = - \sum_{n=1}^{\infty} \frac{\log n }{n^{2}} \cos \left( 
\frac{2 \pi n}{5} \right),
\end{equation}
\noindent
and Theorem \ref{main-thm}  give the identity
\begin{equation}
\ione \frac{\log(1  - \varphi  x + x^{2})}{2x} \lx \, dx = 
\frac{\gamma \pi^{2}}{150} + 
\sum_{n=1}^{\infty} \frac{\log n }{n^{2}} \cos \left( 
\frac{2 \pi n}{5} \right). 
\label{nice-cos}
\st
\end{equation}
\end{example}

\medskip

\begin{example}
\label{bell}
The result of Theorem \ref{main-thm} reduces the evaluation of a certain 
class of integrals to the evaluation of the corresponding $L$-functions. 
Many natural choices of the function $Q$  lead to series that the authors 
are unable to evaluate. For example, $Q(x) = e^{x}$ produces the identity
\begin{equation}
\ione e^{x} \lx \, dx = -\gamma(e-1) - \sum_{n=1}^{\infty} \frac{\log n }{n!},
\st
\end{equation}
\noindent
and we have been unable to procude an analytic expression for the series
above. The same is true for the series in (\ref{nice-cos}).
\end{example}

\section{Evaluation of integrals  with real poles. An expression for
$R_{m,0}(1)$} 
\label{sec-special}
\setcounter{equation}{0}

We now turn to the evaluation of the integrals
\begin{equation}
R_{m,j}(a) = \ione \frac{x^{j} \, \lx }{(x+a)^{m+1}} \, dx.
\end{equation}
\noindent
The method of partial fractions can then be used to produce explicit 
formulas for integrals of the type
\begin{equation}
I_{Q} = \ione Q(x) \, \lx \, dx 
\end{equation}
\noindent
where $Q$ is a rational function with only real poles. \\

In this section we introduce a special family of polynomials $T_{m}(x)$
and produce an explicit analytic expression for
\begin{equation}
E_{m} := \ione \frac{T_{m-1}(x) \, \lx }{(x+1)^{m+1}} \, dx. 
\end{equation}
\noindent
These are then employed to evaluate $R_{m,0}(1)$. \\

\begin{Definition}
The {\em Eulerian polynomials} $A_{m}$ are defined by the generating function
\begin{equation}
\frac{1-x}{1-x \, \text{exp}[t(1-x)]} = \sum_{m=0}^{\infty} A_{m}(x) 
\frac{t^{m}}{m!}.
\label{gen-fun}
\st
\end{equation}
\end{Definition}

\begin{Note}
The Eulerian polynomials appear in many combinatorial problems. The 
coefficients $A_{m,j}$ in
\begin{equation}
A_{m}(x) = \sum_{j=1}^{m} A_{m,j}x^{j}
\end{equation}
\noindent
are the {\em Eulerian numbers}. They count the number of permutations of 
$\{ 1, \, 2, \, \cdots, n \}$ which show exactly $j$ increases between
adjacent elements, the first element always being counted as a jump. The 
numbers $A_{m,j}$ have an explicit formula 
\begin{equation}
A_{m,j} = \sum_{k=0}^{j} (-1)^{k} \binom{m+1}{k} (j-k)^{m},
\label{eulerian-1}
\st
\end{equation}
\noindent
and a recurrence relation
\begin{equation}
A_{m,j} = jA_{m-1,j} + (m-j+1)A_{m-1,j-1}
\st
\end{equation}
\noindent
that follows from 
\begin{equation}
A_{m+1}(x) = x(1-x) \frac{d}{dx} A_{m}(x) + (m+1)x A_{m}(x),
\st
\label{recu-1}
\end{equation}
\noindent
with $A_{0}(x) = 1$. The recurrence (\ref{recu-1}) follows directly from 
(\ref{gen-fun}) and it immediatly implies that $A_{m}(x)$ is a polynomial 
of degree $m$. The first few are
\begin{eqnarray}
A_{0}(x) & = & 1, \nonumber \\
A_{1}(x) & = & x, \nonumber \\
A_{2}(x) & = & x^{2}+ x, \nonumber \\
A_{3}(x) & = & x^{3} + 4x^{2}+ x, \nonumber \\
A_{4}(x) & = & x^{4} + 11x^{3} + 11x^{2}+ x. \nonumber 
\st
\end{eqnarray}
\noindent
More information about these polynomials can be found in \cite{graham1}.
\end{Note}

We now present the relation between Eulerian polynomials and the 
polylogarithm function. 

\begin{Lem}
\label{lem-poly}
Let $m \in \mathbb{N}$. The polynomial $A_{m}(x)$ satisfies
\begin{equation}
\text{PolyLog}[-m,-x] = \frac{A_{m}(-x)}{(x+1)^{m+1}}.
\st
\end{equation}
\end{Lem}
\begin{proof}
The identity 
\begin{equation}
\text{PolyLog}[-m,-x] = \sum_{n=0}^{\infty} (-1)^{n+1} (n+1)^{m} x^{n+1},
\st
\end{equation}
\noindent
shows that 
\begin{equation}
\text{PolyLog}[-m,-x] = \vartheta^{(m)} \left( \frac{1}{1+x} \right),
\st
\end{equation}
\noindent
where $\vartheta = x \frac{d}{dx}$. The claim now follows by using 
(\ref{recu-1}) and an elementary induction.
\end{proof}

We now employ Eulerian polynomials to evaluate an auxiliary 
family of integrals.

\begin{Prop}
Let $m \in \mathbb{N}$. Define
\begin{equation}
T_{m}(x) := - \frac{A_{m+1}(-x)}{x} = \sum_{j=0}^{m} (-1)^{j} A_{m+1,j+1}x^{j},
\label{p-def}
\end{equation}
\noindent
and
\begin{equation}
E_{m} := \ione \frac{T_{m-1}(x) \, \lx }{(x+1)^{m+1}} \, dx. 
\end{equation}
\noindent
Then 
\begin{equation}
E_{m} = (1- 2^{m}) \zeta'(1-m) + \left( \gamma(2^{m}-1) + 2^{m} \log 2 \right) 
\zeta(1-m), \quad \text{ for } m \geq 1.
\st
\end{equation}
\end{Prop}
\begin{proof}
This is a special case of Example \ref{fourth}.
\end{proof}

\begin{Note}
The relation (\ref{zeta-neg}) gives
\begin{equation}
E_{m} = (1-2^{m})\zeta'(1-m) +(-1)^{m+1} \left( \gamma ( 2^{m}-1) + 2^{m} \log 2 \right)
\frac{B_{m}}{m}, \quad \text{ for } m  \geq 1.
\label{expr-em}
\st
\end{equation}
\noindent
For example,
\begin{eqnarray}
E_{1} & = & - \frac{\gamma}{2} + \frac{\log \pi}{2} - \frac{\log 2}{2}, \st \\
E_{2} & = & - \frac{1}{4} -\frac{\gamma}{4} - \frac{\log 2}{3} + 3 \log A, 
\st
\nonumber
\end{eqnarray}
\noindent
where $A$ is the Glaisher constant defined in (\ref{glaisher}).
\end{Note}

The expression for $E_{m}$ in (\ref{expr-em}) is now used to provide a 
recurrence for the integrals $R_{m,0}(1)$, where $R_{m,j}(a)$ is defined in
(\ref{Bmj}). We begin
with an example that will provide an initial condition for the 
recurrence in Corollary \ref{cor-recu}.

\begin{example}
\label{r00}
The integral $R_{0,0}(1)$ is given by
\begin{equation}
R_{0,0}(1) = \ione \frac{\lx \, dx}{1+x} = - \frac{\log^{2} 2}{2}.
\st
\end{equation}
\noindent
The choice $a_{n} = (-1)^n$ in Theorem \ref{main-thm} 
produces $Q(x) = 1/(1+x)$ and $L_{Q}(s) = (1- 2^{1-s})\zeta(s)$. Passing to
the limit as $s \to 1$ and using $R_{0,0}(1) = - \gamma L_{Q}(1) + 
L_{Q}'(1)$ gives the result. 
\end{example}

\begin{Thm}
The integrals $E_{m}$ in (\ref{expr-em}) satisfy
\begin{equation}
E_{m} = \sum_{k=0}^{m-1} b_{k}(m) R_{m-k,0}(1),
\st
\end{equation}
\noindent
where 
\begin{equation}
b_{k}(m) = (-1)^{k} \sum_{j=k}^{m-1} \binom{j}{k} A_{m,j+1}
\end{equation}
\noindent
and $A_{m,j}$ are the Eulerian numbers given in (\ref{eulerian-1}). 
\end{Thm}
\begin{proof}
In the expression
\begin{equation}
E_{m} = \ione \frac{T_{m-1}(x) \, \lx }{(x+1)^{m+1}} \, dx
\end{equation}
\noindent
use (\ref{p-def}) to obtain 
\begin{equation}
E_{m} = \sum_{j=0}^{m-1} (-1)^{j} A_{m,j+1} \ione \frac{x^{j} \, \lx }
{(x+1)^{m+1}} \, dx.
\end{equation}
\noindent
Now write $x = (x+1)-1$, expand the resulting binomial and reverse the order 
of summation to obtain the result. 
\end{proof}

\begin{Cor}
\label{cor-recu}
The integrals $R_{m,0}(1)$ satisfy the recurrence 
\begin{eqnarray}
R_{1,0}(1) & = & E_{1} \label{initial-1} \st \\
R_{m,0}(1) & = & \frac{E_{m}}{b_{0}(m)} - 
\sum_{k=1}^{m-1} \frac{b_{k}(m)}{b_{0}(m)} R_{m-k,0}(1).
\label{initial-2}
\st
\end{eqnarray}
\end{Cor}
\begin{proof}
First observe that $b_{0}(m) \neq 0$. Indeed, 
\begin{equation}
b_{0}(m) = \sum_{j=1}^{m} A_{m,j} = A_{m}(1).
\st
\end{equation}
Using the recurrrence (\ref{recu-1}) we conclude that 
$A_{m+1}(1) = (m+1)A_{m}(1)$.  Therefore
$b_{0}(m) = A_{m}(1) = m!$. 
\end{proof}

\begin{example}
\label{r01-example}
The previous result provides the values 
\begin{eqnarray}
R_{1,0}(1) & = & \frac{1}{2} \left( -\gamma + \log \pi - \log 2 \right), 
\st
\nonumber \\
R_{2,0}(1) & = & \frac{1}{24} 
\left( -3 - 9 \gamma - 10 \log 2 + 36 \log A + 6 \log \pi \right) 
\st
\nonumber \\
R_{3,0}(1) & = & 
\frac{1}{24} \left( -3 - 7 \gamma - 8 \log 2 + 36 \log A + 4 \log \pi + 
7 \zeta(3)/\pi^{2} \right).
\st
\nonumber
\end{eqnarray}
\noindent
These integrals are computable using Mathematica 6.0. 
\end{example}

\section{An expression for $R_{m,0}(a)$} \label{sec-bmzeroa}
\setcounter{equation}{0}

In this section we present an analytic expression for 
\begin{equation}
R_{m,0}(a) := \ione \frac{\lx }{(x+a)^{m+1}} \, dx.
\end{equation}
\noindent
The result is given in terms of the polylogarithm function defined in 
(\ref{polylog-def0})  and the derivative 
$\text{PolyLog}^{(1,0)}[c,x]$ defined in (\ref{polyder-def0}).
The evaluation employs 
the expression for $R_{m,0}(1)$ given the previous section.

\begin{Prop}
\label{prop-51}
The integral $R_{0,0}(a)$ is given by 
\begin{equation}
R_{0,0}(a) := \ione \frac{\lx \, dx}{x+a} = - \gamma \log (1 + 1/a) - 
\text{PolyLog}^{(1,0)}[1,-1/a].
\st
\end{equation}
\end{Prop}
\begin{proof}
The expansion
\begin{equation}
Q(x) = \frac{1}{x+a} = \frac{1}{a} \sum_{n=0}^{\infty} (-1)^{n} a^{-n} x^{n},
\st
\end{equation}
\noindent
produces the $L$-function
\begin{equation}
L_{Q}(s) = 
\sum_{n=0}^{\infty} \frac{ (-1)^{n}}{a^{n+1} (n+1)^{s}} = 
- \text{PolyLog}[s, - 1/a ].
\st
\end{equation}
\noindent
Theorem \ref{main-thm} gives the result.
\end{proof}

The evaluation of $R_{m,0}(a)$ for $m \geq 1$ employs the {\em signless 
Stirling numbers of the first kind} $S_{1}(m,j)$ defined by the expansion
\begin{equation}
(t)_{m} = \sum_{j=1}^{m} S_{1}(m,j) t^{j},
\label{stir-1}
\end{equation}
\noindent
where $(t)_{m} = t(t+1)(t+2) \cdots (t+m-1)$ is the Pochhammer symbol. 

\begin{Thm}
Let $m \in \mathbb{N}$ and $a>0$. The 
$L$-function associated to $Q(x) = 1/(x+a)^{m+1}$
is 
\begin{equation}
L_{Q}(s) = - \frac{1}{a^{m} m!} \sum_{j=1}^{m} S_{1}(m,j) 
\text{PolyLog}[s-j, - 1/a ].
\end{equation}
\end{Thm}
\begin{proof}
The identity
\begin{equation}
\binom{-\beta}{k} = \frac{(-1)^{k} (\beta)_{k}}{k!}
\st
\end{equation}
\noindent
is used in the expansion
\begin{equation}
\frac{1}{(x+a)^{\beta}} = a^{-\beta} ( 1 + x/a)^{-\beta} = 
a^{-\beta} \sum_{k=0}^{\infty} \binom{-\beta}{k} a^{-k} x^{k}
\st
\end{equation}
\noindent
to produce
\begin{equation}
\frac{1}{(x+ a)^{\beta}} = a^{-\beta} \sum_{k=0}^{\infty} 
\frac{(- 1/a)^{k} (\beta)_{k}}{k! (k+1)^{s}}. 
\end{equation}
\noindent
Now choose $\beta = m+1$ and use the elementary identity
\begin{equation}
\frac{(m+1)_{k}}{k!} = \frac{(k+1)_{m}}{m!},
\st
\end{equation}
\noindent
to write the $L$-function corresponding to $Q(x) = 1/(x+a)^{m+1}$ as 
\begin{equation}
L_{Q}(s) = \frac{a^{-(m+1)}}{m!} 
\sum_{k=0}^{\infty} \frac{(-1/a)^{k} (k+1)_{m}}{(k+1)^{s}}.
\end{equation}
\noindent
Finally use the expression (\ref{stir-1}) to write
\begin{equation}
L_{Q}(s) = \frac{a^{-(m+1)}}{m!} 
\sum_{j=1}^{m} S_{1}(m,j) \sum_{k=0}^{\infty} 
\frac{( - 1/a)^{k}}{(k+1)^{s-j}},
\end{equation}
\noindent
and identify the series as a polylogarithm to produce the result.
\end{proof}

\begin{Cor}
\label{cor-53}
Let $m \in \mathbb{N}$. Then the integral
\begin{equation}
R_{m,0}(a) := \ione \frac{\lx }{(x+a)^{m+1}} \, dx 
\end{equation}
\noindent
is given by
\begin{multline}
R_{m,0}(a)  =  - \frac{\gamma}{a^{m}(1+a)m} - 
\frac{\gamma}{a^{m+1} m!} \sum_{j=2}^{m} \frac{S_{1}(m,j) T_{j-2}(1/a)}
{(1+ 1/a)^{j}} \nonumber \\
  - \frac{1}{a^{m} m!} 
\sum_{j=1}^{m} S_{1}(m,j) \text{PolyLog}^{(1,0)}[1-j,-1/a],
\st
\nonumber
\end{multline}
\noindent
where $T_{j}$ is the polynomial defined in (\ref{p-def}).
\end{Cor}
\begin{proof}
The result follows from Theorem \ref{main-thm} and the expression for 
polylogarithms given in (\ref{lem-poly}).
\end{proof}

\begin{example}
The choice $a=2$ and $m=1$ gives
\begin{equation}
R_{1,0}(2) := \ione \frac{\lx \, dx}{(x+2)^{2}} = - \frac{\gamma}{6} 
- \frac{1}{2} \text{PolyLog}^{(1,0)}[0, -\tfrac{1}{2}].
\st
\end{equation}
\noindent
Mathematica 6.0 is unable to compute these integral.  The expansion of the 
polylogarithm function gives the identity 
\begin{equation}
\ione \frac{\lx \, dx}{(x+2)^{2}} = - \frac{\gamma}{6} + 
\sum_{n=1}^{\infty} \frac{(-1)^{n} \, \log n }{2^{n+1}}.
\st
\end{equation}
\end{example}

\section{An algorithm for the evaluation of $R_{m,j}(a)$} \label{sec-bmja}
\setcounter{equation}{0}

In this section we present an analytic expression for 
\begin{equation}
R_{m,j}(a) := \ione \frac{x^{j} \, \lx }{(x+a)^{m+1}}  \, dx.
\end{equation}
\noindent
These results provide the evaluation of integrals of the type
\begin{equation}
I_{Q} = \ione Q(x) \lx \, dx
\end{equation}
\noindent
where $Q$ has real poles. The case of non-real poles is discussed in the 
next section.

\begin{Thm}
\label{bmj-rec}
The integrals $R_{m,j}(a)$ satisfy the recurrence
\begin{equation}
R_{m,0}(a) = \sum_{j=0}^{r} \alpha_{j,r}(a) R_{m-r+j,j}(a),
\label{triangular}
\end{equation}
\noindent
for any $r \leq m$. Here
\begin{equation}
\alpha_{j,r}(a) := (-1)^{j} \binom{r}{j} a^{-r}.
\label{alpha}
\end{equation}
\end{Thm}

\begin{Note}
The recurrence (\ref{triangular}) can now be used for 
increasing values of the free parameter $r$, to obtain 
analytic expressions for $R_{m,j}(a)$. For instance, $r=1$ gives 
\begin{equation}
R_{m,0}(a) = \alpha_{0,1}(a)R_{m-1,0}(a) + \alpha_{1,1}(a) R_{m,1}(a), 
\end{equation}
\noindent
that determines $R_{m,1}(a)$ in terms of $R_{m,0}(a)$ and $R_{m-1,0}(a)$, that
were previously computed. 
\end{Note}

The proof of Theorem \ref{bmj-rec} employs a recurrence for the functions 
$\alpha_{j,r}(a)$ that is established first.

\begin{Lem}
\label{lemma-0}
Let $k, r \in \mathbb{N}$ and $\alpha_{j,r}(a)$ as in (\ref{alpha}). Then 
\begin{equation}
\frac{\alpha_{0,r}(a)}{(x+a)^{k}} + 
\frac{\alpha_{1,r}(a) x}{(x+a)^{k+1}} +  \cdots + 
\frac{\alpha_{r,r}(a) x^{r}}{(x+a)^{k+r}}  = 
\frac{1}{(x+a)^{k+r}}. 
\end{equation}
\end{Lem}
\begin{proof}
Expand the identity
\begin{equation}
1 = \frac{(x+a)^{r}}{a^{r}} \left( 1 - \frac{x}{x+a} \right)^{r}.
\nonumber
\end{equation}
\end{proof}

\noindent
{\em Proof of Theorem \ref{bmj-rec}}. Multiply the relation in Lemma 
\ref{lemma-0} by $\lx$ and integrate over $[0,1]$. \\

\begin{example}
We now use the method described above to check that
\begin{equation}
R_{1,1}(1) := \ione \frac{x \, \lx }{(x+1)^{2}} \, dx =
\frac{1}{2} \left( -\log^{2}2 + \gamma - \log \pi + \log 2 \right).
\st
\end{equation}
\noindent
The recursion (\ref{triangular}) gives 
\begin{equation}
R_{1,0}(1) = \alpha_{0,1}(1) R_{0,0}(1) + \alpha_{1,1}(1) R_{1,1}(1).
\st
\end{equation}
\noindent
Using the values $\alpha_{0,1}(1) = 1$ and 
$\alpha_{1,1}(1) = -1$ and the integrals 
\begin{equation}
R_{0,0}(1) = - \frac{1}{2} \log^{2}2 
\st
\end{equation}
\noindent
given in Example \ref{r00} and 
\begin{equation}
R_{1,0}(1) = \frac{1}{2} \left( - \gamma + \log \pi - \log 2 \right),
\st
\end{equation}
\noindent
computed in Example \ref{r01-example}, we obtain the  result.
\end{example}

\begin{example}
The computation of $R_{1,1}(2)$ can be obtained from the recurrence
\begin{equation}
R_{1,0}(2) = \alpha_{0,1}(2) R_{0,0}(2) + \alpha_{1,1}(2)R_{1,1}(2),
\st
\end{equation}
\noindent
and the previously computed values 
\begin{equation}
R_{1,0}(2) = - \frac{\gamma}{6} - 
\frac{1}{2} \text{PolyLog}^{(1,0)}[0, - \tfrac{1}{2} ],
\st
\end{equation}
\noindent
and 
\begin{equation}
R_{0,0}(2) = - \gamma \log \frac{3}{2}  -
\text{PolyLog}^{(1,0)}[1, - \tfrac{1}{2} ].
\st
\end{equation}
\noindent
It follows that
\begin{eqnarray}
R_{1,1}(2) & := & \int_{0}^{1} \frac{x \, \lx }{(x+2)^{2}} \, dx \nonumber \\
& = & \frac{\gamma}{3} - \gamma \log \frac{3}{2} +
\text{PolyLog}^{(1,0)}[0, - \tfrac{1}{2} ] 
- \text{PolyLog}^{(1,0)}[1, - \tfrac{1}{2} ] \st \nonumber \\
& = & \frac{\gamma}{3} - \gamma \log \frac{3}{2} -
\sum_{n=2}^{\infty} \frac{(-1)^{n} \log n}{2^{n}} ( 1 - 1/n ).
\st
\nonumber
\end{eqnarray}
\end{example}

\begin{example}
We now illustrate the recurrence (\ref{triangular}) to obtain the value 
\begin{equation}
R_{3,2}(5)  =   \ione \frac{x^{2} \, \lx }{(x+5)^{4}} \, dx.
\end{equation}

We first let $m=3$ in (\ref{form-132}) to obtain
\begin{equation}
R_{3,0}(5) = \alpha_{0,2}(5) R_{1,0}(5) + 
\alpha_{1,2}(5) R_{2,1}(5) + 
\alpha_{2,2}(5) R_{3,2}(5).
\label{mess-1}
\st
\end{equation}
\noindent
The integrals with second index $0$ are given in (\ref{cor-53}) by
\begin{equation}
R_{1,0}(5) = - \frac{\gamma}{30} - \frac{1}{5}\text{PolyLog}^{(1,0)}[0, 
- \tfrac{1}{5} ] 
\label{rzero1-5}
\st
\end{equation}
\noindent
and 
\begin{multline}
R_{3,0}(5)  =   - \frac{91 \gamma}{81000} - 
\frac{1}{750}\text{PolyLog}^{(1,0)}[-2, - \tfrac{1}{5} ]  \nonumber \\
  - \frac{1}{250}\text{PolyLog}^{(1,0)}[-1, - \tfrac{1}{5} ] -
\frac{1}{375}\text{PolyLog}^{(1,0)}[0, - \tfrac{1}{5} ].
\st
\nonumber 
\end{multline}
\noindent
The next step is to put $m=2$ in (\ref{form-131}) to obtain
\begin{equation}
R_{2,0}(5) = \alpha_{0,1}(5) R_{1,0}(5) + \alpha_{1,1}(5) R_{2,1}(5).
\label{mess-2}
\st
\end{equation}
\noindent
The values 
\begin{equation}
R_{2,0}(5) = - \frac{11 \gamma}{1800} 
-\frac{1}{50}\text{PolyLog}^{(1,0)}[-1, - \tfrac{1}{5} ]
-\frac{1}{50}\text{PolyLog}^{(1,0)}[0, - \tfrac{1}{5} ]
\nonumber
\st
\end{equation}
\noindent
and $R_{1,0}(5)$ is given in (\ref{rzero1-5}). These come from 
Corollary \ref{cor-53}. Equation (\ref{mess-2}) now gives 
\begin{equation}
R_{2,1}(5) = - \frac{\gamma}{360} +
\frac{1}{10}\text{PolyLog}^{(1,0)}[-1, - \tfrac{1}{5} ] -
\frac{1}{10}\text{PolyLog}^{(1,0)}[0, - \tfrac{1}{5} ].
\nonumber
\st
\end{equation}
\noindent
Finally we obtain 
\begin{multline}
R_{3,2}(5)  =  -\frac{\gamma}{3240} 
-\frac{1}{30}\text{PolyLog}^{(1,0)}[-2, - \tfrac{1}{5} ] \nonumber \\
 +  \frac{1}{10}\text{PolyLog}^{(1,0)}[-1, - \tfrac{1}{5}  ]
-\frac{1}{15}\text{PolyLog}^{(1,0)}[0, - \tfrac{1}{5} ] 
\st
\nonumber 
\end{multline}
\noindent
from (\ref{mess-1}). 
\end{example}

\medskip

\begin{Note}
The integrals $R_{m,j}(a)$ are computable by Mathematica 6.0 for $a=1$, but 
not for $a \neq 1$. 
\end{Note}

\medskip

\section{Evaluation of integrals with non-real poles. The 
integrals $C_{m,j}(a,b) = D_{m,j}(r, \theta)$} \label{sec-complex}
\setcounter{equation}{0}

We consider now the evaluation of integrals 
\begin{equation}
C_{m,j}(a,b) := \ione \frac{x^{j}}{(x^{2} + ax+b)^{m+1}}
 \lx  \, dx, \label{Cmj-1}
\end{equation}
\noindent
where $a^{2}-4b < 0$, so that the quadratic factor has non-real zeros. This is 
written as 
\begin{equation}
x^{2} + ax + b = (x-c)(x- \bar{c}) = x^{2} - 2rx \cos \theta  + r^{2}
\end{equation}
\noindent
and we write
\begin{equation}
\label{def-D}
D_{m,j}(r, \theta) = \ione \frac{x^{j}}{(x^{2}-2rx \cos \theta  + r^{2})^{m+1}}
\, \lx \, dx. 
\end{equation}
\noindent
Naturally $C_{m,j}(a,b) = D_{m,j}(r, \theta)$, we are simply emphasizing the 
polar representation of the poles.  \\

\noindent
{\bf Plan of evaluation}: the computation of $D_{m,j}(r, \theta)$ can be 
reduced to the range $m \geq 0$ and $0 \leq j \leq 2m+1$ by dividing $x^{j}$ 
by $(x^{2}-2rx \cos \theta + 1)^{m+1}$, in case $j \geq 2m+2$. The fact is that 
the recurrences (\ref{recu-one}) and (\ref{recu-two}) determine all the 
integrals $D_{m,j}(r, \theta)$ from 
$D_{0,0}(r, \theta)$ and $D_{0,1}(r, \theta)$. This is illustrated with the 
four integrals $D_{1,j}(r, \theta): \, 0 \leq j \leq 3$. Begin with 
(\ref{recu-one})
with $m=j=1$. This  gives 
\begin{equation}
D_{1,1}(r, \theta) = - \frac{1}{2r \sin \theta} 
\frac{\partial}{\partial \theta} D_{0,0}(r, \theta)
\end{equation}
\noindent
and then (\ref{recu-two}) with $m = j = 1$ gives
\begin{equation}
D_{1,1}(r, \theta) = \frac{1}{2 \cos \theta} 
\left( \frac{\partial}{\partial r} D_{0,0}(r, \theta) + 2 r D_{1,0}(r, \theta)
\right),
\end{equation}
\noindent
and this determines $D_{1,0}(r, \theta)$.  Now use $m=1,\, j=2$ in 
(\ref{recu-one}) to obtain
\begin{equation}
D_{1,2}(r, \theta) = - \frac{1}{2r \sin \theta} 
\frac{\partial}{\partial \theta} D_{0,1}(r, \theta).
\end{equation}
\noindent
Finally, (\ref{recu-two}) with $m=1$ and $j=3$ yields
\begin{equation}
D_{1,3}(r, \theta) = \frac{1}{2 \cos \theta} 
\left( \frac{\partial }{\partial r } D_{0,2}(r, \theta)
 + 2r D_{1,2}(r, \theta) \right),
\end{equation}
\noindent
Dividing $x^{2}$ by $x^{2} - 2r x \cos \theta + r^{2}$ expresses 
$D_{0,2}(r, \theta)$ as a linear combination of $D_{0,0}(r, \theta)$ and 
$D_{1,0}(r, \theta)$. This process determines $D_{1,3}(r, \theta)$.  \\

We compute first the integral $D_{0,0}(r, \theta)$. This task is 
is divided into two cases according to whether $r=1$ or not.
Theorem \ref{thm-73} gives the result for the first case and 
Theorem \ref{thm-78} describes the case $r \neq 1$. 
The evaluation of the integrals $D_{m,j}(r, \theta)$ are then obtained by 
using the recurrences (\ref{recu-one}) and (\ref{recu-two}). \\

\noindent
{\bf Calculation of $D_{0,0}(1, \theta)$}. This is stated in the next 
theorem.  \\

\begin{Thm}
\label{thm-73}
Assume $0 < \theta < 2 \pi$. Then 
\begin{eqnarray}
D_{0,0}(1, \theta ) & := & 
\ione \frac{\lx }{x^{2} - 2x \cos \theta + 1} \, dx \nonumber \\
& = & \frac{\pi}{2 \sin \theta} 
\left[ ( 1 - \theta/\pi) \log 2 \pi + 
\log \left( \frac{\Gamma(1- \theta/ 2 \pi)}{\Gamma( \theta/ 2 \pi)} 
\right) \right].
\nonumber 
\st
\end{eqnarray}
\end{Thm}
\begin{proof}
Consider the function
\begin{equation}
Q(x) = \frac{\sin \theta}{x^{2} - 2x \cos \theta + 1} 
\end{equation}
\noindent
with the classical expansion
\begin{equation}
Q(x) = \sum_{k=0}^{\infty} \sin \left( (k+1) \theta \right) x^{k}.
\end{equation}
\noindent
The corresponding $L$-function is given by
\begin{equation}
L_{Q}(s) = \sum_{k=0}^{\infty} \frac{\sin \left( (k+1) \theta \right)}
{(k+1)^{s}}, 
\end{equation}
\noindent
and its value at $s=1$ is given by
\begin{equation}
L_{Q}(1) = 
\sum_{k=0}^{\infty} \frac{\sin \left( (k+1) \theta \right)}{k+1} =
\frac{\pi - \theta}{2},
\end{equation}
\noindent
while the closed form of the derivative at $s=1$ is 
\begin{eqnarray}
L_{Q}'(1) & = & - \sum_{k=0}^{\infty} \frac{\sin \left( (k+1) \theta 
\right) }{k+1} \log(k+1) \label{der-1} \\
& = & - \frac{\pi}{2} 
\left( \log \left( \frac{\Gamma(\theta/ 2 \pi)}{\Gamma( 1 - \theta/ 2 \pi)}
\right) + ( \gamma + \log 2 \pi ) \left( \frac{\theta}{\pi} -1 \right) \right).
\nonumber
\end{eqnarray}
\noindent
This identity can be found in \cite{apelblat}, page 250, $\# 30$.  The result
now follows from Theorem  \ref{main-thm}.
\end{proof}

\medskip

\noindent
{\bf Calculation of $D_{0,0}(r, \theta)$ in the case $r \neq 1$}. This is 
stated in the theorem below.  \\

\begin{Thm}
\label{thm-78}
Assume $0 < \theta < 2 \pi$ and $r \neq 1$. Then
\begin{eqnarray}
D_{0,0}(r, \theta)  & := &  \ione \frac{ \lx }{x^{2} - 2rx \cos \theta + r^{2}} 
\, dx \nonumber \\
&  =  & - \frac{\gamma}{r \sin \theta} 
\tan^{-1} \left( \frac{\sin \theta}{r - \cos \theta}  \right) \nonumber \\
   & & \quad  +  \, \frac{1}{2 r i \sin \theta} 
\left( 
\text{PolyLog}^{(1,0)}[1, e^{i \theta}/r] - 
\text{PolyLog}^{(1,0)}[1, e^{-i \theta}/r] 
\right). 
\nonumber 
\end{eqnarray}
\end{Thm}
\begin{proof}
Consider the function
\begin{equation}
Q(x) = \frac{r^{2} \sin \theta}{x^{2} - 2rx \cos \theta + r^{2}} 
= \sum_{k=0}^{\infty} \frac{\sin \left( (k+1) \theta \right) }{r^{k}} x^{k}
\end{equation}
\noindent
and the corresponding $L$-function
\begin{eqnarray}
L_{Q}(s) & := & \sum_{k=0}^{\infty} \frac{\sin \left( (k+1) \theta \right) }
{r^{k} (k+1)^{s}} \nonumber \\
& = & \sum_{k=0}^{\infty} \frac{e^{i(k+1) \theta} - e^{-i(k+1) \theta}}
{2 i r^{k} (k+1)^{s}} \nonumber \\
& = & \frac{r}{2i} \left( \text{PolyLog}[s, e^{it}/r] - 
\text{PolyLog}[s, e^{-it}/r]  \right). \nonumber 
\end{eqnarray}
\noindent
The identity 
\begin{equation}
\text{PolyLog}[1,a] = - \log(1-a),
\end{equation}
\noindent
yields the value
\begin{equation}
L_{Q}(1) = \frac{r}{2i} \left[ - \log( 1 - e^{i \theta}/r ) + 
\log(1 - e^{-i \theta}/r ) \right]. 
\label{l-complex}
\end{equation}
\noindent
The value of $L_{Q}(1)$ can be written as 
\begin{equation}
L_{Q}(1) = r \tan^{-1} \left( \frac{\sin \theta}{r - \cos \theta} \right). 
\label{l-real}
\end{equation}
\noindent
This follows by checking that both expressions for $L_{Q}(1)$ 
match at $\theta = 0$ and their derivatives with respect to $\theta$ match. 
\end{proof}

We now present several special cases of these evaluations. Many of them 
appear in the table of integrals \cite{gr}. \\

\begin{example}
The value $D_{0,0}(1, \theta)$ in Theorem \ref{thm-73} 
appears as $4.325.7$ in \cite{gr}. 
\end{example}

\begin{example}
Replacing $\theta$ by $\theta + \pi$, we obtain the evaluation
\begin{equation}
\ione \frac{\lx }{x^{2} + 2x \cos \theta + 1} \, dx =
 \frac{\pi}{2 \sin \theta} 
\left[ \frac{\theta \log 2 \pi}{\pi} + 
\log \left( \frac{\Gamma(1/2 + \theta/2 \pi)}{\Gamma(1/2 - \theta/2 \pi)} 
\right) \right].
\nonumber 
\st
\end{equation}
\noindent
This appears as $4.231.2$ in \cite{gr}.
\end{example}

\begin{example}
The value $\theta = \pi/2$ provides 
\begin{equation}
D_{0,0} \left(1 , \frac{\pi}{2} \right) = \ione \frac{\lx }{1+x^{2}} \, dx = \frac{\pi}{4} \log 2 \pi + 
\frac{\pi}{2} \log \frac{ \Gamma(3/4) }{\Gamma(1/4)}.
\st
\nonumber
\end{equation}
\noindent
This is the example discussed by Vardi in \cite{vardi1}.
\end{example}

\begin{example}
The angle $\theta = \pi/3$ yields the value 
\begin{equation}
D_{0,0} \left( 1, \frac{\pi}{3} \right) = \ione \frac{\lx}{1-x+x^{2}} \, dx = 
\frac{2 \pi \log 2 \pi}{3 \sqrt{3}} + \frac{\pi}{\sqrt{3}} 
\log \frac{ \Gamma(5/6) }{\Gamma(1/6)}.
\nonumber
\st
\end{equation}
\noindent
This appears as $4.325.6$ in \cite{gr} where the answer is written in the 
equivalent form
\begin{equation}
\ione \frac{\lx}{1-x+x^{2}} \, dx = 
\frac{2 \pi }{\sqrt{3}} \left[  \frac{5}{6} \log 2 \pi - 
\log \Gamma(1/6) \right].
\st
\nonumber
\end{equation}
\noindent
The equivalent form
\begin{equation}
\ift \frac{\log x \, dx}{e^{x} + e^{-x} -1} = 
\frac{2 \pi }{\sqrt{3}} \left[  \frac{5}{6} \log 2 \pi - 
\log \Gamma(1/6) \right]
\st
\nonumber
\end{equation}
\noindent
appears as $4.332.1$ in \cite{gr}. 
\end{example}

\begin{example}
The angle $\theta = 2 \pi/3$ provides an evaluation of $4.325.5$ in 
\cite{gr}:
\begin{equation}
D_{0,0} \left( 1, \frac{ 2 \pi}{3} \right) = 
\ione \frac{\lx}{1+x+x^{2}} \, dx = 
\frac{\pi \log 2 \pi}{3 \sqrt{3}} + \frac{\pi}{\sqrt{3}} 
\log \frac{ \Gamma(2/3) }{\Gamma(1/3)}.
\st
\nonumber
\end{equation}
\noindent
The equivalent form 
\begin{equation}
\ift \frac{\log x \, dx}{e^{x} + e^{-x} +1} = 
\frac{\pi}{\sqrt{3}} \log \left( \frac{\Gamma(2/3) \, \sqrt{2 \pi}}
{\Gamma(1/3)} \right)
\end{equation}
\noindent
appears {\em incorrectly} as $4.332.2$ in \cite{gr}. The correct result is
obtained by replacing $\sqrt{2 \pi}$  by $\sqrt[3]{2 \pi}$. 
\end{example}

\begin{example}
The limit of $D_{0,0}(1, \theta)$ as $\theta \to \pi$ gives the evaluation of
\begin{equation}
D_{0,0} \left( 1, \frac{\pi}{2} \right) = 
\ione \frac{\lx }{(1+x)^{2}} \, dx = \frac{1}{2}
\left( \log \pi - \log 2 - \gamma \right). 
\st
\nonumber
\end{equation}
\noindent
This is $4.325.3$ of \cite{gr}.
\end{example}

\begin{example}
It is easy to choose an angle and produce an integral that cannot be 
evaluated by Mathematica 6.0. For example, $\theta = 3 \pi /4$ gives 
\begin{equation}
D_{0,0} \left(1, \frac{3 \pi}{4} \right) = 
\ione \frac{\lx}{1 + \sqrt{2} \, x + x^{2}} \, dx = 
\frac{\pi}{\sqrt{2}} \left( \frac{\log 2 \pi}{4} + 
\log \left( \frac{\Gamma(5/8)}{\Gamma(3/8)} \right) \right).
\st
\nonumber
\end{equation}
\end{example}

\medskip

\begin{Note}
The evaluation of $L_{Q}(1)$ yields the identity
\begin{equation}
\label{nice-0}
\ione \frac{dx}{x^{2} - 2rx \cos \theta + r^{2}} = 
\frac{1}{r \sin \theta} 
\tan^{-1} \left( \frac{\sin \theta}{r - \cos \theta} 
\right).
\st
\end{equation}
\end{Note}

\begin{Note}
The polylogarithm terms appearing in the expression for $D_{0,0}(r, \theta)$
can be written as in terms of the sum
\begin{equation}
U(r, \theta) := \sum_{k=0}^{\infty} \sin \left( (k+1) \theta \right) 
\, \frac{\log(k+1)}{r^{k} (k+1)}.
\end{equation}
\noindent
The authors are unable to express the function $U(r, \theta)$ in terms of 
special functions with real arguments. 
\end{Note}

We now proceed to a systematic determination of the integrals 
$D_{m,j}(r, \theta)$  for $m, \, j > 0$. For that, we use the recurrences 
(\ref{recu-one}) and (\ref{recu-two}).

\begin{Prop}
Assume $0 < \theta < 2 \pi$. Then
\begin{eqnarray}
D_{1,1}(1, \theta ) & := & 
\ione \frac{x \, \lx }{(x^{2} - 2 x \cos \theta  +1 )^{2}} \, dx  \nonumber \\
& = & \frac{\log 2 \pi}{4 \sin^{2} \theta} ( 1 + ( \pi - \theta) \cot \theta )
+ \frac{1}{8 \sin^{2} \theta} \left( \psi \left( \theta/ 2 \pi \right) + 
\psi \left( 1 - \theta/ 2 \pi \right) \right) \nonumber \\
& & + \, \frac{\pi}{4} \csc^{2} \theta \cot \theta 
\log \left( \frac{\Gamma( 1 - \theta/2 \pi) }{\Gamma( \theta/ 2 \pi)} \right). 
\nonumber 
\end{eqnarray}
\end{Prop}
\begin{proof}
The result follows directly from 
(\ref{recu-one}) and Theorem \ref{thm-73}.
\end{proof}

\begin{Note}
For the case $r \neq 1$ the value of $D_{1,1}(r, \theta)$ can be obtained 
by differentiating the expression for $D_{0,0}(r, \theta)$ in Theorem 
\ref{thm-78}.
\end{Note}

\medskip

Particular cases of this result are stated next. \\

\begin{example}
\label{example-77}
The angle $\theta = \pi/2$ produces
\begin{equation}
D_{1,1} \left(1, \frac{\pi}{2} \right) := 
\ione \frac{x \, \lx }{(x^2+1)^{2}} \, dx = \frac{\log 2 \pi}{4} + 
\frac{1}{8} \psi \left( \frac{1}{4} \right) + 
\frac{1}{8} \psi \left(\frac{3}{4} \right). 
\st
\nonumber
\end{equation}
\noindent
Here $\psi(x) = \Gamma'(x)/\Gamma(x)$ is the polygamma function. Using the 
values 
\begin{equation}
\psi \left( \frac{1}{4} \right) = - \gamma - \frac{\pi}{2} - 3 \log 2 
\text{ and }
\psi \left( \frac{3}{4} \right) = - \gamma + \frac{\pi}{2} - 3 \log 2 
\nonumber 
\st
\end{equation}
\noindent
that appear in \cite{gr} as $8.366.4$ and $8.366.5$ respectively, we obtain 
\begin{equation}
\ione \frac{x \, \lx }{(x^2+1)^{2}} \, dx = 
\frac{1}{4} ( \log \pi - 2 \log 2 - \gamma).
\st
\end{equation}
\noindent
The equivalent version
\begin{equation}
\ift \frac{\log x \, dx }{\cosh^{2}x} = 
\log \pi - 2 \log 2 - \gamma,
\st
\end{equation}
\noindent
appears as $4.371.3$ in \cite{gr}. 
\end{example}

\begin{example}
The angle $\theta = \pi/3$ gives the evaluation
\begin{eqnarray}
D_{1,1} \left(1, \frac{\pi}{3} \right)& := & 
\ione \frac{x \, \lx }{(x^{2} - x + 1)^{2}} \, dx  \st \nonumber \\
& = & \frac{\log 2 \pi}{3} + \frac{2 \pi \log 2 \pi}{9 \sqrt{3}} + 
\frac{\pi}{3 \sqrt{3}} \log \left( \frac{ \Gamma( 5/6) }{\Gamma(1/6)} \right)
\nonumber \\ 
&  & + \, \frac{1}{6} \psi \left(\frac{1}{6} \right) + 
\frac{1}{6} \psi \left(\frac{5}{6} \right). 
\nonumber 
\end{eqnarray}
\noindent
Using elementary properties of the $\psi$ function and the values 
\begin{equation}
\psi \left( \frac{1}{6} \right) = - \gamma - \frac{\pi \sqrt{3} }{2} - 
\frac{3}{2} \log 3 - 2 \log 2,  
\st
\end{equation}
\noindent
and 
\begin{equation}
\psi \left( \frac{5}{6} \right) =  - \gamma + \frac{\pi \sqrt{3} }{2} - 
\frac{3}{2} \log 3 - 2 \log 2, 
\st
\end{equation}
\noindent
that appear in \cite{srichoi}, page $21$, we obtain
\begin{eqnarray}
\ione \frac{x \, \lx }{(x^{2} - x + 1)^{2}} \, dx & = & 
- \frac{\gamma}{3} - \frac{\log 2}{3} + \frac{7 \pi \log 2}{9 \sqrt{3}} - 
\frac{\log 3}{2} - \frac{\pi \log 3}{3 \sqrt{3}} \nonumber \\
 &  & + \frac{\log \pi}{3} + \frac{8 \pi \log \pi}{9 \sqrt{3}} - 
\frac{4 \pi}{3 \sqrt{3}} \log \Gamma \left( \frac{1}{3} \right).
\nonumber
\st
\end{eqnarray}

\end{example}

\begin{example}
The recurrences (\ref{recu-one}) and (\ref{recu-two}) yield
\begin{eqnarray}
\ione \frac{x^{2} \, \lx }{(x^{2} - \sqrt{3}x + 1)^{3}} \, dx & = & 
\frac{1}{6}( 9 \sqrt{3} + 25 \pi) \log 2 \pi + 
5 \pi \log \left( \frac{ \Gamma(11/12)}{\Gamma(1/12)} \right) \st \nonumber \\
& + & \frac{3 \sqrt{3}}{4} \left( \psi \left(1/12 \right) + 
\psi \left(11/12 \right) \right) + 
\frac{1}{8 \pi} \left( \psi' \left(11/12 \right) - 
\psi' \left(1/12 \right)  \right). \nonumber
\end{eqnarray}
\noindent
This can be written as 
\begin{eqnarray}
\ione \frac{x^{2} \, \lx }{(x^{2} - \sqrt{3}x + 1)^{3}} \, dx & = & 
- \frac{3 \sqrt{3} \gamma}{2} - 3 \sqrt{3} \log 2 + \frac{35}{3} \pi \log 2 
 - \frac{9}{4} \sqrt{3} \log 3 + \frac{9}{2} \log ( 2 -  \sqrt{3} ) 
\nonumber \\
& & - 5 \pi \log(\sqrt{3} -1 ) + \frac{3 \sqrt{3}}{2} \log \pi 
+ \frac{55}{6} \pi \log \pi - 10 \pi \log \Gamma \left(1/12 \right) 
\nonumber \\
& &  
- \frac{1}{8 \pi} \psi' \left(1/12 \right) 
+  \frac{1}{8 \pi} \psi'\left( 11/12 \right).
\st
\nonumber 
\end{eqnarray}
\end{example}

\begin{example}
The values $r=2$ and $\theta = \pi/3$ yields the evaluation
\begin{equation}
\ione \frac{\lx \, dx}{x^{2} - 2x + 4}  =  - \frac{\gamma \pi}{6 \sqrt{3}} 
-  \frac{i}{ 2 \sqrt{3}} 
\left(
\text{PolyLog}^{(1,0)} \left[ 1 , \tfrac{1+ i \sqrt{3}}{4} \right] - 
\text{PolyLog}^{(1,0)} \left[ 1 , \tfrac{1- i \sqrt{3}}{4} \right] 
\right).
\nonumber
\end{equation}
\noindent
Mathematica 6.0 is unable to evaluate this integral.
\end{example}

\noindent
{\bf Calculation of $D_{0,1}(r, \theta)$}.  \\

The integral 
\begin{equation}
D_{0,1}(r, \theta) = \ione \frac{x \lx \, dx}{x^{2}-2rx \cos \theta + r^{2}}
\end{equation}
\noindent
corresponds to 
\begin{equation}
Q(x)  = \frac{x}{x^{2}-2rx \cos \theta + r^{2}}.
\end{equation}
\noindent 
To evaluate the integral $D_{0,1}(r, \theta)$ we employ the expansion 
\begin{equation}
\frac{x}{x^{2} - 2 r x \cos \theta + r^{2}} = \sum_{k=0}^{\infty} 
\frac{\sin k \theta}{r^{k+1} \, \sin \theta} x^{k}. 
\label{key-1}
\end{equation}
\noindent
We conclude that the $L$-function associated to this $Q$ is 
\begin{equation}
L_{Q}(s) = \frac{1}{\sin \theta} 
\sum_{k=0}^{\infty} \frac{\sin k \theta }{r^{k+1} 
\, (k+1)^{s}}. 
\label{L-sum-0}
\end{equation}
\noindent
Therefore 
\begin{equation}
L_{Q}(1) = \frac{1}{\sin \theta} 
\sum_{k=0}^{\infty} \frac{\sin (k+1) \theta}{r^{k+1} \, 
(k+2)}.
\end{equation}
\noindent
To evaluate this sum, integrate (\ref{key-1}) from $0$ to $1$ to produce 
\begin{equation}
L_{Q}(1) = \ione \frac{x \, dx}{x^{2}-2rx \cos \theta + r^{2}}.
\end{equation}
\noindent
Observe that
\begin{eqnarray}
\ione \frac{x \, dx}{x^{2}-2rx \cos \theta + r^{2}} & = & 
\frac{1}{2} \ione \frac{(2x - 2r \cos \theta) \, dx}
{x^{2}-2rx \cos \theta + r^{2}}  +  \nonumber \\
 & + & r \cos \theta 
\ione \frac{dx}{x^{2}-2rx \cos \theta + r^{2}}. \nonumber 
\end{eqnarray}
\noindent
Both integrals are elementary, the latter is given in (\ref{nice-0}). Therefore,
\begin{equation}
L_{Q}(1) = \frac{1}{2} \log \left( \frac{r^{2}- 2r  \cos \theta +1}{r^{2}} 
\right) + \cot \theta \, \tan^{-1} \left( \frac{\sin \theta }{r - \cos \theta}
\right). 
\end{equation}

The $L$-series (\ref{L-sum-0}) can be expressed in terms of the  Lerch 
$\Phi$-function 
\begin{equation}
\Phi(z,s,a) := \sum_{k=0}^{\infty} \frac{z^{k}}{(k+a)^{s}}.
\end{equation}
\noindent
Indeed, 
\begin{eqnarray}
L_{Q}(s) & = & \frac{1}{2i \sin \theta} \sum_{k=0}^{\infty} 
\frac{e^{ik \theta} - e^{-i k \theta}}{r^{k+1} \, (k+1)^{s}} 
\nonumber \\
& = & \frac{1}{2i r \sin \theta} 
\left( \sum_{k=0}^{\infty} \frac{( e^{i \theta}/r)^{k}}{(k+1)^{s}} -
\frac{( e^{- i \theta}/r)^{k}}{(k+1)^{s}} \right) \nonumber \\
& = & \frac{1}{2ri \sin \theta } 
\left[ \Phi \left( \frac{e^{i \theta}}{r}, s , 1 \right) - 
       \Phi \left( \frac{e^{- i \theta}}{r}, s , 1 \right) \right]. 
\nonumber
\end{eqnarray}

The next statement gives the value of $D_{0,1}(r, \theta)$.

\begin{Thm}
The integral
\begin{equation}
D_{0,1}(r, \theta) = \ione \frac{x \lx \, dx}{x^{2}-2xr \cos \theta + r^{2}}
\end{equation}
\noindent
is given by
\begin{eqnarray}
D_{0,1}(r, \theta) & = & - \frac{\gamma}{2} \log \left( \frac{r^{2} - 
2 r \cos \theta +1}{r^{2}} \right) \nonumber \\
& - & \gamma \cot \theta \, \tan^{-1} \left( \frac{\sin \theta}
{r - \cos \theta } \right) \nonumber \\
& + & 
\frac{1}{2ri \sin \theta } 
\left[ \Phi^{(0,1,0)} \left( \frac{e^{i \theta}}{r}, 1 , 1 \right) - 
       \Phi^{(0,1,0)} \left( \frac{e^{- i \theta}}{r}, 1 , 1 \right) \right]. 
\nonumber
\end{eqnarray}
\end{Thm}

\medskip

\begin{Note}
This evaluation completes the algorithm to evaluate all the integrals 
$D_{m,j}(r, \theta)$. 
\end{Note}

\section{Adamchik's integrals} \label{sec-adamchik}
\setcounter{equation}{0}

V. Adamchik presented in \cite{adamchik3} a series of beautiful evaluations 
of integrals of the form (\ref{type-1}), where the denominator has the form
$(1+x^{n})^{m}$ for $n, \, m \in \mathbb{N}$. The results are expressed in 
terms of the Hurwitz zeta function and its derivatives. For example, in 
Proposition $3$, it is shown that
\begin{eqnarray}
\ione \frac{x^{p-1}}{1+x^{n}} \, \lx \, dx & = & 
\frac{\gamma + \log 2n }{2n} 
\left( \psi \left( \frac{p}{2n} \right) - 
       \psi \left( \frac{n+ p}{2n} \right) \right) \nonumber \\
& + & \frac{1}{2n} \left( \zeta' \left(1, \frac{p}{2n} \right) -  
                   \zeta' \left(1, \frac{n+p}{2n} \right) \right),
\nonumber 
\end{eqnarray}
\noindent
followed by Proposition $4$ that states that
\begin{eqnarray}
\ione x^{p-1}\frac{1-x}{1-x^{n}} \, \lx \, dx & = & 
\frac{\gamma + \log n }{n} 
\left( \psi \left( \frac{p}{n} \right) - 
       \psi \left( \frac{p+ 1}{n} \right) \right) \nonumber \\
& + & \frac{1}{n} \left( \zeta' \left(1, \frac{p}{n} \right) -  
                   \zeta' \left(1, \frac{p+1}{n} \right) \right). 
\nonumber 
\end{eqnarray}
\noindent
The expressions  become more complicated as the exponent of the denominator 
increases. For instance, Proposition $5$ gives
\begin{eqnarray}
\ione \frac{x^{p-1}}{(1+x^{n})^{2}} \lx \, dx & = & 
\frac{(n-p)}{2n^{2}} ( \log 2n + \gamma ) 
\left( \psi \left( \frac{p}{2n} \right) - \psi \left( \frac{n+p}{2n} \right) \right) \nonumber \\
& - & \frac{1}{2n} \left( \gamma + \log 2n - 
2 \log \left( \frac{\Gamma(\tfrac{p}{2n}) }{\Gamma(\tfrac{n+p}{2n})} \right) 
\right)
\nonumber \\
& + & \frac{n-p}{2n^{2}} 
\left( \zeta' \left( 1, \frac{p}{2n} \right) - 
        \zeta' \left( 1, \frac{n+p}{2n} \right) \right), \nonumber
\end{eqnarray}
\noindent
and in Proposition $6$ we find
\begin{eqnarray}
\ione \frac{x^{p-1}}{(1+x^{n})^{3}} \, \lx \, dx & = & 
\frac{3n-2p}{2n^{2}} \log \left( \frac{\Gamma\left( \frac{p}{2n} \right)}
                                      {\Gamma\left( \frac{n+p}{2n} \right)}
\right) 
- \frac{(5n-2p)( \log 2n + \gamma)}{8n^{2}}  \nonumber \\
& + & \frac{(n-p)(2n-p)( \log 2n + \gamma)}{4n^{2}} 
\left( \psi \left( \frac{p}{2n} \right) - \psi \left( \frac{n+p}{2n} \right) \right) \nonumber \\
& + & \frac{1}{n} 
\left( \zeta' \left( -1, \frac{p}{2n} \right) - 
        \zeta' \left( -1, \frac{n+p}{2n} \right) \right) \nonumber \\
& + & \frac{(n-p)(2n-p)}{4n^{3}} 
\left( \zeta' \left( 1, \frac{p}{2n} \right) - 
        \zeta' \left( 1, \frac{n+p}{2n} \right) \right). \nonumber
\end{eqnarray}

We now describe some examples on how to use 
Theorem \ref{main-thm} to obtain some of the 
specific examples in \cite{adamchik3}. 

\begin{example}
From Proposition $3$ it follows that
\begin{equation}
\ione \frac{x^{n-1}}{1+x^{n}} \, \lx \, dx = 
- \frac{\log 2 \, \log 2n^{2}}{2n}.
\label{formula-27}
\end{equation}
\noindent
This appears as $(27)$ in \cite{adamchik3}. To check this evaluation, observe 
that 
\begin{equation}
Q(x) = \frac{x^{n-1}}{1+ x^{n}} = \sum_{k=1}^{\infty} (-1)^{k-1} x^{kn-1},
\end{equation}
\noindent
so the corresponding $L$-function is 
\begin{equation}
L_{Q}(s) = \sum_{k=1}^{\infty} \frac{(-1)^{k-1}}{k^{s}n^{s}} = 
(1 - 2^{-s}) \frac{\zeta(s)}{n^{s}}.
\end{equation}
A direct calculation shows that $L_{Q}(1) = \frac{\log 2}{n}$ and 
\begin{equation}
L_{Q}'(1) = \frac{\gamma \log 2}{n} - \frac{\log^{2}2}{2n} - 
\frac{\log 2 \, \log n}{n}. 
\end{equation}
Then (\ref{formula-27}) follows from (\ref{main-formula}).
\end{example}

\begin{example}
Formula $(28)$ in \cite{adamchik3} is also obtained from Proposition $3$ 
and it states that
\begin{equation}
\ione \frac{x^{2n-1}}{1+x^{n}} \, \lx \, dx = 
\frac{1}{2n} \left( \log^{2}2 + 2( \log 2 -1) \log n - 2 \gamma \right).
\end{equation}
\noindent
To prove this, consider 
\begin{equation}
Q(x) = \frac{x^{2n-1}}{1+x^{n}} = \sum_{k=0}^{\infty} (-1)^{k} x^{(2+k)n-1},
\end{equation}
\noindent
whose $L$-function is 
\begin{equation}
L_{Q}(s) = \sum_{k=0}^{\infty} \frac{(-1)^{k}}{(2+k)^{s} n^{s}} = 
\frac{1}{n^{s}} \left( 1 - ( 1 - 2^{1-s}) \zeta(s) \right).
\end{equation}
\noindent
Replacing the  values
\begin{equation}
L_{Q}(1) = \frac{1 - \log 2 }{n} \text{ and }
L_{Q}'(1) = \frac{1}{2n} \left( \log^{2} 2 - 2 \gamma \log 2 - 2 \log n + 
2 \log 2 \log n \right),
\nonumber
\end{equation}
\noindent
in (\ref{main-formula}) we obtain the result.
\end{example}

\begin{example}
The identity
\begin{equation}
\ione \frac{x}{1+x^{4}} \, \lx \, dx = 
\frac{\pi}{4} \log \left( \frac{ \sqrt{\pi} \, \Gamma( \tfrac{3}{4} ) }
{\Gamma( \tfrac{1}{4} ) } \right),
\label{last-adam}
\end{equation}
\noindent
appears as formula $(30)$ in \cite{adamchik3}. To establish it,
consider the function
\begin{equation}
Q(x) = \frac{x}{1+x^{4}} = \sum_{k=0}^{\infty} (-1)^{k} x^{4k+1},
\end{equation}
\noindent
with $L$-function
\begin{equation}
L_{Q}(s) = \sum_{k=0}^{\infty} \frac{(-1)^{k}}{(4k+2)^{s}} = 
\frac{1}{2^{3s}} \left( \zeta(s, \tfrac{1}{4}) - \zeta(s, \tfrac{3}{4}) \right).
\st
\end{equation}
\noindent
The result follows from Theorem \ref{main-thm} by using the values
\begin{equation}
L_{Q}(1) = \tfrac{1}{8} \left( \psi \left( \tfrac{3}{4} \right) - 
                              \psi \left( \tfrac{1}{4} \right) \right), 
\end{equation}
\noindent
and 
\begin{equation}
L_{Q}'(1) = - \tfrac{3 \log 2}{8} 
\left( \psi \left( \tfrac{3}{4} \right) - 
               \psi \left( \tfrac{1}{4} \right) \right) + 
\tfrac{1}{8} \left( \zeta'(1, \tfrac{1}{4}) - 
                   \zeta'(1, \tfrac{3}{4}) \right). 
\end{equation}
\noindent
The special values 
\begin{equation}
\psi(\tfrac{1}{4} ) =  -\gamma - \tfrac{\pi}{2} - 3 \log 2 \text{ and }
\psi(\tfrac{3}{4} ) =  -\gamma + \tfrac{\pi}{2} - 3 \log 2, 
\end{equation}
\noindent
and the value  
\begin{equation}
\zeta'(1,\tfrac{1}{4} ) - 
\zeta'(1,\tfrac{3}{4} ) = 
\pi \left( \gamma + \log 2 + 3 \log 2 \pi - 4 \log \Gamma \left( \tfrac{1}{4}
\right) \right), 
\end{equation}
\noindent
are used to simplified the result. This 
last expression comes from
\begin{eqnarray}
\zeta'(1, \tfrac{p}{q} ) - \zeta'(1, 1 - \tfrac{p}{q} )  & =  & 
\pi \cot \frac{\pi p}{q} \left( \log 2 \pi q + \gamma \right) \nonumber \\
& - & 2 \pi \sum_{j=1}^{q-1} \log \left( \Gamma \left( \tfrac{j}{q} \right) 
\right) \sin \frac{2 \pi j p}{q}.
\nonumber
\end{eqnarray}
\noindent
This identity, established in \cite{adamchik3}, follows directly from the 
classical Rademacher formula
\begin{equation}
\zeta \left(z, \tfrac{p}{q} \right) = 2 \Gamma(1-z) (2 \pi q)^{z-1} 
\sum_{j=1}^{q} \sin \left( \tfrac{\pi z}{2} + \tfrac{2jp z}{q} \right) 
\zeta \left( 1 -z , \tfrac{j}{q} \right).
\end{equation}

\medskip

An alternative evaluation of this integral comes from the partial fraction 
decomposition 
\begin{equation}
\frac{x}{x^{4} + (2 - a^{2})x^{2} +  1} = 
\frac{1}{2a} \frac{1}{x^{2}-ax+1} -
\frac{1}{2a} \frac{1}{x^{2}+ax+1}. 
\label{pfd}
\end{equation}
\noindent 
We assume $|a| < 2 $ and write $a = 2 \cos \theta$. Then (\ref{pfd}) 
yields 
\begin{equation}
\ione \frac{x \, \lx}{x^{4} + (2 - 4 \cos^{2} \theta)x^{2} + 1} \, dx =
\frac{1}{4 \cos \theta} \left( D_{0,0}(1, \theta) - D_{0,0}(1, \pi + \theta)
\right). 
\end{equation}
\noindent 
Using the result of Theorem \ref{thm-73} we obtain
\begin{multline}
\ione \frac{x \, \lx \, dx}{x^{4} + (2 - 4 \cos^{2} \theta)x^{2} + 1}  = 
\frac{\pi}{4 \sin 2 \theta}  \times \label{nice} \st  \\ 
\left( 
\log \left( \frac{4 \pi^{3}}{\sin \theta} \right) - \frac{2 \theta}{\pi} 
\log 2 \pi - 
2 \log  \left[ \Gamma \left( \frac{\theta}{2 \pi} \right) 
       \Gamma \left( \frac{1}{2} + \frac{\theta}{2 \pi} \right) \right] 
\right). 
\end{multline}
\noindent
The  special case $\theta = \pi/4$ produces (\ref{last-adam}).
\end{example}

\begin{example}
The case $\theta = \pi/2$ in the previous example reduces 
to Example \ref{example-77}. 
\end{example}

\begin{example}
The angle  $\theta = \pi/3$ in (\ref{nice}) yields
\begin{equation}
\ione \frac{x \, \lx }{x^{4} + x^{2} + 1 } \, dx = 
\frac{\pi}{12 \sqrt{3}} 
\left( 6 \log 2 - 3 \log 3 + 8 \log \pi - 12 \log \Gamma \left( \tfrac{1}{3} 
\right) \right).
\end{equation}
\end{example}

\begin{example}
The angle  $\theta = \pi/8$ in (\ref{nice}) yields
\begin{equation}
\ione \frac{x \, \lx }{x^{4}  - \sqrt{2}  x^{2} + 1 } \, dx = 
\frac{\pi}{8  \sqrt{2}} 
\left( 7 \log \pi - 4 \log \sin \tfrac{\pi}{8} 
-8  \log \Gamma \left( \tfrac{1}{8} 
\right) \right).
\end{equation}
\end{example}

\section{A hyperbolic example} \label{sec-hyper}
\setcounter{equation}{0}

The method introduced here can be used to provide an analytic expression
for the family 
\begin{equation}
LC_{n}:= \ift \frac{\log t \, dt}{\cosh^{n+1}t}.
\end{equation}
\noindent
The table of integrals \cite{gr} contains
\begin{equation}
LC_{0}= \ift \frac{\log t \, dt}{\cosh t} = 
\frac{\pi}{2} \left( 2 \log 2 + 3 \log \pi - 4 \log \Gamma \left( 
\tfrac{1}{4} \right) \right),
\st
\end{equation}
\noindent
as formula $4.371.1$ and
\begin{equation}
LC_{1}= \ift \frac{\log t \, dt}{\cosh^{2}t} =  
- \gamma + \log \pi - 2 \log 2,
\st
\end{equation}
\noindent
as $4.371.3$.  \\

The change of variables $x = e^{-t}$ shows that
\begin{equation}
LC_{n} = 2^{n+1} \ione \frac{x^{n} \, \lx }{(x^{2} + 1)^{n+1}} \, dx,
\end{equation}
\noindent
that identifies this integral as
\begin{equation}
LC_{n} = 2^{n+1} D_{n,n} \left(1, \tfrac{\pi}{2} \right)
\end{equation}

The recurrence (\ref{recu-one}), for $r=1$ and $j=m$, become
\begin{equation}
D_{m,m}(1, \theta)  =  - \frac{1}{2m \sin \theta} 
\frac{\partial}{\partial \theta}
D_{m-1,m-1}(1, \theta), 
\end{equation}
\noindent
and the initial condition 
\begin{equation}
D_{0,0}(1, \theta ) 
 =  \frac{\pi}{2 \sin \theta} 
\left[ ( 1 - \theta/\pi) \log 2 \pi + 
\log \left( \frac{\Gamma(1- \theta/ 2 \pi)}{\Gamma( \theta/ 2 \pi)} 
\right) \right],
\nonumber 
\end{equation}
\noindent
provides a systematic procedure to compute $LC_{n}$. For instance, it follows
that
\begin{eqnarray}
LC_{2} & = & 2^{3} D_{2,2}(1, \tfrac{\pi}{2})  \nonumber \\ 
     & = & - \frac{2}{\pi} \text{Catalan} + 
\frac{\pi}{4} \left( 2 \log 2 + 3 \log \pi - 4 \log \Gamma \left( 
\tfrac{1}{4} \right) \right),
\st
\nonumber
\end{eqnarray}
\noindent
and 
\begin{eqnarray}
LC_{3} & = & 2^{4} D_{3,3}(1, \tfrac{\pi}{2})  \nonumber \\ 
     & = & - \frac{2 \gamma }{3} - \frac{ 4 \log 2}{3} + 
\frac{2 \log \pi}{3} + \frac{28}{3} \zeta'(-2). \nonumber
\st
\end{eqnarray}
\noindent
The {\em Catalan constant} appearing above is defined by 
\begin{equation}
\text{Catalan} = \sum_{n=0}^{\infty} \frac{(-1)^{n}}{(2n+1)^{2}}.
\st
\end{equation}

\section{Small sample of a new type of evaluations} \label{sec-sample}
\setcounter{equation}{0}

We have introduced here a systematic method to deal with integrals of the form
\begin{equation}
I_{Q} = \ione Q(x) \lx \, dx.
\end{equation}
\noindent
Extensions of this technique provides examples such as 
\begin{eqnarray}
\ione (1 + \log x) \log(x+1)  \, \lx \, dx & = & 
1 + (\gamma-1) \left(\frac{\pi^2}{12} - 1 \right)  \st \nonumber \\
& & \, - \left( \frac{\pi^{2}}{12} + 2 \right) \log 2   
 - \, \frac{\zeta'(2)}{2}, \nonumber \\
\ione ( 1 + \log x) \, \tan^{-1}x  \, \lx \, dx & = & 
(1 - \gamma) \frac{\pi^{2}}{48} - \frac{\pi}{4} + \frac{\log 2}{2} + 
\frac{\zeta'(2)}{8}, \st \nonumber \\
\ione \frac{\tanh^{-1}\sqrt{x}}{x} \, \lx \, dx & = & 
- \frac{\gamma \pi^{2}}{4} + \frac{\pi^{2} \log 2}{3} + \frac{3}{2} \zeta'(2).
\st
\nonumber 
\end{eqnarray}

\medskip

\noindent
Details will presented elsewhere.

\section{Conclusions} \label{sec-conclusions}
\setcounter{equation}{0}

We have developed an algorithm to evaluate integrals of the form
\begin{equation}
I_{Q} = \ione Q(x) \, \lx \, dx.
\end{equation}
\noindent
In the case where $Q(x)$ is analytic at $x=0$, with power series expansion
\begin{equation}
Q(x) = \sum_{n=0}^{\infty} a_{n} x^{n}
\end{equation}
\noindent
we associate to $Q$ its $L$-function 
\begin{equation}
L_{Q}(s) = \sum_{n=0}^{\infty} \frac{a_{n}}{(n+1)^{s}}.
\end{equation}
\noindent
Then 
\begin{equation}
I_{Q} = - \gamma L_{Q}(1) + L_{Q}'(1).
\end{equation}

In the case $Q(x)$ is a rational function, we provide explicit expressions for 
$I_{Q}$ in terms of special values of the logarithm, the Riemann zeta function, 
the polylogarithm $\text{PolyLog}[c,x]$, its first derivative with respect 
to $c$ and the Lerch $\Phi$-function. 

\bigskip

\noindent
{\bf Acknowledgments}. The second
author was partially funded by
$\text{NSF-DMS } 0070567$. The first author was supported by 
VIGRE grant $\text{NSF-DMS } 0239996$ as a graduate student.

\end{document}